\documentclass[10pt]{article}

\usepackage{amssymb}
\usepackage{amsfonts}
\usepackage{amsmath}
\usepackage{enumerate}
\usepackage{theorem}
\usepackage{ifthen}
\usepackage[scaled=0.92]{helvet}
\usepackage[latin1]{inputenc}

\newcommand{\R}{{\mathbb {R}}}

\newcommand{\1}{{\bf 1}}

\newcommand{\buil}[3]{\mathrel{\mathop{\kern0pt#1}\limits_{#2}^{#3}}}

\def\tresp#1_#2^#3{\mathrel {\mathop{\kern 0pt#1}\limits_{#2}^{#3}}}

\theoremheaderfont{\bfseries}
\newtheorem{theorem}{Theorem}

\newtheorem{proposition}[theorem]{Proposition}
\newtheorem{lemma}[theorem]{Lemma}

{\theorembodyfont{\rmfamily}}

\renewcommand{\theequation}{\thesection.\arabic{equation}}

\def\tresp#1_#2^#3{\mathrel {\mathop{\kern 0pt#1}\limits_{#2}^{#3}}}

\def\^#1{\if#1i{\accent"5E\i}\else{\accent"5E #1}\fi}
\def\"#1{\if#1i{\accent"7F\i}\else{\accent"7F #1}\fi}
\def\1{{\rm l}\hskip -0.21truecm 1}

\begin{document}
\thispagestyle{empty}
\null
\begin{center}
{\LARGE\bf Anticipating Linear Stochastic Differential Equations Driven by a L\'{e}vy Process}
\end{center}
\vspace{1cm}

\noindent {\large\bf Jorge A. Le\'{o}n}
 
\noindent Departamento de Control Automático, Cinvestav-IPN

\noindent Apartado Postal 14-740, 07000 M\'{e}xico D.F., Mexico

\noindent jleon@ctrl.cinvestav.mx
\vspace{0.5cm}

\noindent {\large \bf David M\'{a}rquez-Carreras}

\noindent Departament de Probabilitat, L\`{o}gica i Estad\'{\i}stica, Facultat de Matem\`atiques, Universitat de Barcelona

\noindent Gran Via 585, 08007-Barcelona (Catalunya)

\noindent davidmarquez@ub.edu
\vspace{0.5cm}

\noindent {\large \bf Josep Vives}

\noindent Departament de Probabilitat, L\`{o}gica i Estad\'{\i}stica, Facultat de Matem\`atiques, Universitat de Barcelona

\noindent Gran Via 585, 08007-Barcelona (Catalunya), Spain

\noindent josep.vives@ub.edu
\vspace{1cm}


\begin{abstract}
\noindent
In this paper we study the existence of a unique solution for  linear stochastic differential equations driven by a L\'evy process, where the initial condition and the coefficients are random and not necessarily adapted to the underlying filtration. Towards this end, we  extend the method based on Girsanov
transformation on Wiener space and developped by Buckdahn
\cite{bu3} to the canonical L\'evy space, which is introduced
in \cite{SUV07}.

\vspace{0.3cm}

\vspace{0.3cm}


\end{abstract}

\newpage

\renewcommand{\theequation}{\thesection.\arabic{equation}}

\section{Introduction}
\setcounter{equation}{0}
\setcounter{theorem}{0}

Our  aim in this paper is to prove the existence and uniqueness of a solution of the  linear stochastic differential equation 
\begin{eqnarray}\label{eintr2}
X_t&=&X_0+ \int_0^t b_s X_s\ ds +  \int_0^t a_s X_s\ \delta W_s + \int_0^t \int_{\{|y|>1\}} v_{s-}(y) X_{s-}\ dN(s,y)\nonumber\\
&&+\int_0^t \int_{\{0<|y|\leq 1\}} v_{s-}(y) X_{s-}\ d\tilde N(s,y),\qquad 0\le t \le T.
\end{eqnarray}
Here $X_0$ is a random variable, $a$, $b$ and $v(y)$, for any $y\in{\mathbb R}$, $y\neq 0,$ are random processes no necessarily adapted to the underlying filtration, $W$ is the canonical Wiener process, $N$ is the canonical Poisson random measure with parameter $\nu$ (see Section \ref{sec:clp} for details), $d{\tilde N}(t,y):=dN(t,y)-dt\,\nu(dy)$, and the integral with respect to $W$ (respectively the integrals with respect to $N$ and ${\tilde N}$) is in the Skorohod sense (respectively are pathwise defined). 

In the adapted case (i.e., deterministic initial condition and adapted coefficients to the filtration generated by $W$ and $N$), the stochastic differential equation (\ref{eintr2}) with
 no necessarily linear coefficients has been analyzed by several authors  (see, for instance, \cite{Ap,BL,BW,FO,IW,Ku,Ja,Pr,Ru}). For example, Ikeda and Watanabe \cite{IW} have considered 
this equation with no necessarily linear coefficients and have used the Picard iteration procedure and Gronwall´s lemma to show existence and uniqueness of the solution, respectively. It is well-known that this is possible due to the isometry property of It\^o integrals. Also, in this case, one approach to study equation (\ref{eintr2}) is to assume first that $N$ does not have small jumps (i.e., the absolutely values of the jumps side are bigger than a constant $\varepsilon>0$) and consider equation (\ref{eintr2}) as an stochastic differential equation driven by a Brownian motion between two consecutive jump times, which has a unique solution under suitable conditions due to It\^o \cite{It}. Then, we only need to show that this solution converges to the one of equation (\ref{eintr2}) as $\varepsilon\rightarrow 0$. Namely, the solution of the equation
\begin{eqnarray}\label{eintr3}
X_t^\varepsilon&=&X_0+ \int_0^t b_s X_s^\varepsilon\ ds + \int_0^t a_s X_s^\varepsilon\ \delta W_s + \int_0^t \int_{\{|y|>1\}} v_{s-}(y) X_{s-}^\varepsilon\ dN(s,y)\nonumber\\
&&+\int_0^t \int_{\{\varepsilon<|y|\leq 1\}} v_{s-}(y) X_{s-}^\varepsilon\ d\tilde N(s,y),\qquad 0\le t \le T,
\end{eqnarray}
converges as $\varepsilon\downarrow 0$, to a solution of equation (\ref{eintr2}). We can see Rubenthaler \cite{Ru} for details. This method was also utilized to obtain an It\^o formula for L\'evy processes (see, for example, Cont and Tankov \cite{CT04}). We also mention that in the adapted and linear case, It\^o formula provides a tool to obtain the existence and uniqueness of  the solution to (\ref{eintr2}). For details, the reader can consult Protter \cite{Pr}.

In the general case, we cannot use neither the Picard iteration procedure, nor Gronwall´s lemma to deal with (\ref{eintr2}) because the $L^2$-norm of the solution depends on its derivative in the Malliavin calculus sense and this derivatives can be estimated only in terms of the second derivative, and so on. Therefore we do not have a closed argument, as it is pointed out by Nualart \cite{N06}.
 
On the Wiener case (i.e., $v\equiv 0$), Buckdahn \cite{Bu,Bu2,bu3} has study equation (\ref{eintr2}) via anticipating Girsanov transformations. In particular, he showed that It\^o formula is not useful in this case. This approach has been also useful to deal
with fractional stochastic differential equations (see
\cite{JM,JJ}).

On the Poisson space, it means $a\equiv 0$, equation (\ref{eintr2}) has been considered in different situations for different definitions of stochastic integral (see, for instance, \cite{Jjj,JT,LRT,Jpi,Npr}).

In this paper, in order to  obtain the existence of a unique solution to equation (\ref{eintr2}), we  apply the method developed in 
\cite{Bu,Bu2,bu3} between consecutive jumps times to figure out the
solution of the  stochastic linear equation (\ref{eintr3}). Then, we
get the convergence of $X^\varepsilon$ to the solution of  (\ref{eintr2}).

The paper is organized as follows. Section 2 is devoted to different preliminaires: Canonical Lévy space and process, Malliavin calculus and  anticipative Girsanov transformations. In section 3 the solution candidates for equations (\ref{eintr2}) and 
(\ref{eintr3}) are presented and some of their properties are pointed out. In section 4 the existence of a unique solution of 
(\ref{eintr3}) is proved and in Section 5, the same is done for (\ref{eintr2}). A long and non-central proof of Theorem \ref{continuity} is placed in the Appendix. 

\section{Preliminaries}
\setcounter{equation}{0}
\setcounter{theorem}{0}

 In this section we give the framework and the tools we use in this paper to study the existence of a unique solution to equation (\ref{eintr2}). In particular we introduce the canonical L\'evy space as it was done in Sol\'e et al. \cite{SUV07}, we extend some results given in Buckdahn \cite{Bu,Bu2} to the last space and recall some basic facts of the Malliavin calculus.  
 
In the remaining of this paper, $\nu$ represents a L\'evy measure on $\mathbb R$ such that $\nu(\{0\})=0$ and $\int_{\mathbb R}x^2\nu(dx)<\infty$ (for details see Sato \cite{S99}), $T$ is a positive fixed number and $\ell$ denotes the Lebesgue measure on $[0,T].$ The Borel $\sigma$-algebra of a set $A\subset {\mathbb R}$ is denoted by ${\cal B}(A)$. The jumps of a c\`adl\`ag process $Z$ are denoted by $\Delta Z$ (i.e., $\Delta_t Z=Z_t-Z_{t-}$). Also, for any $p\geq 1$, $|\cdot|_p$ and $||\cdot||_p$ denote the norms on $L^p ([0,T])$ and on $L^p (\Omega)$, respectively.  In particular $||\cdot||_{\infty}$ we will denote the norm on $L^{\infty} (\Omega),$ that is, the essential supremum in $\Omega.$  Sometimes we use the notation $|\cdot|_{L^p(\Omega)}=||\cdot||_p$.

\subsection{Canonical L\'evy space}\label{sec:cls}

 In this paper we consider all the processes defined on the canonical L\'evy space on $[0,T]$,
$$(\Omega,{\cal F},P)=(\Omega_W\otimes \Omega_N, {\cal F}_W\otimes{\cal F}_N,P_W\otimes P_N).$$
Here $(\Omega_W, {\cal F}_W,P_W)$ is the canonical Wiener space and $(\Omega_N, {\cal F}_N P_N)$ is the canonical L\'evy space for a pure jump L\'evy process with L\'evy measure $\nu$, which is defined as follows:

Let $\{\varepsilon_n: n\in{\mathbb N}\}$ be a strictly decreasing sequence of positive numbers such that $\varepsilon_1=1$,
$\lim_{n\rightarrow\infty}\varepsilon_n=0$ and $\nu(S_n)>0$ for any $n\geq 1$, where $S_1=\{x\in{\mathbb R}: \varepsilon_1<|x|\}$ and $S_n=\{x\in{\mathbb R}: \varepsilon_n<|x|\le \varepsilon_{n-1}\}.$ With this notation in mind, the canonical L\'evy space with measure $\nu$ is 
$$(\Omega_N, {\cal F}_N, P_N)=\bigotimes_{n\ge 1}(\Omega^{(n)},{\cal F}^{(n)},P^{(n)}),$$
where $(\Omega^{(n)},{\cal F}^{(n)},P^{(n)})$ is the canonical L\'evy space for a compound Poisson process with intensity $\lambda_n:=\nu(S_n)$ and probability measure $Q_n:=\nu(\cdot\cap S_n)/\lambda_n$. That is, for $n\in{\mathbb N},$
$$\Omega^{(n)}:=\bigcup_{k\ge0}\left([0,T]\times S_n\right)^k,$$
with $\left([0,T]\times S_n\right)^0=\{\alpha\}$, where $\alpha$ is an arbitrary point,
$${\cal F}^{(n)}:=\left\{B\subset\Omega^{(n)}: B\cap\left([0,T]\times S_n\right)^k\in {\cal B}\left(\left([0,T]\times S_n\right)^k\right), \ \hbox{\rm for all } k\in{\mathbb N}\right\}$$
and for any $B\in {\cal F}^{(n)},$
$$P^{(n)}(B):=e^{-\lambda_n T}\sum_{k=0}^{\infty} \frac{\lambda_n^k (\ell \otimes Q_n)^{\otimes k}(B\cap ([0,T]\times S_n)^k)}{k!}.$$

\subsection{Canonical L\'evy process}\label{sec:clp}

 The canonical Wiener process $W=\{W_t:t\in[0,T]\}$ is defined as $W_t(\omega)=\omega(t)$ for 
$\omega \in \Omega_W$, that is, $\omega$ is a continuous function on $[0,T]$ such that $\omega(0)=0$. 

The canonical pure jump process $J_t=\{J_t:t\in[0,T]\}$, with L\'evy measure $\nu$, is
$$J_t(\omega)=\lim_{k\rightarrow\infty}\sum_{n=2}^k\left(X_t^{(n)}(\omega^{(n)})-t\int_{S_n}x\nu(dx)\right)
+X_t^{(1)}(\omega^{(1)}),\quad \omega=(\omega^{(n)})_{n\ge1}\in\Omega_N,$$
where the limit exists with probability 1 and 
$$X_t^{(n)}(\omega^{(n)})=
\begin{cases}
\sum_{l=1}^mx_l1_{[0,t]}(t_l), & \mbox{\rm if $\omega^{(n)}=((t_1,x_1),
\ldots,(t_m,x_m))$,}\\
0, & \mbox{\rm if $\omega^{(n)}=\alpha$}.
\end{cases}
$$
Finally, the canonical Lévy process with triplet $(\gamma, \sigma, \nu)$ is defined as 
$$X_t(\omega)=\gamma t+\sigma W_t(\omega')+J_t(\omega''),\quad\hbox{\rm for } \omega=(\omega',\omega'')\in\Omega_W\otimes\Omega_N.$$ 
Recall also that the associated Poisson random measure is 
$$N(B):=\#\{t\in[0,T]:(t,\Delta X_t)\in B\},\quad B\in{\cal B}([0,T]\times{\mathbb R}_0),$$
where ${\mathbb R}_0={\mathbb R}-\{0\}$.

\subsection{Elements of Malliavin calculus}

 In this paper we deal with the derivative with respect to the process $W$ in the Malliavin calculus sense. So,
in this subsection, we recall some  basic properties of this operator. For details, the reader can consult Nualart\cite{N06} or Sol\'e et al. \cite{SUV07}. 

Let $\mathcal{S}^W$ be the set of random variables of the form 
\begin{equation}\label{eq:rsmo}
F=f(\int_0^T h_1(s)dW_s,\dots,\int_0^T h_n (s)dW_s),
\end{equation}
where $n\geq 1$, $h_j\in L^2([0,T])$ and $f\in \mathcal{C}_{b}^{\infty }(\mathbb{R}^{n}),$ that means  $f$ and all its partial derivatives are bounded. The derivative of the random variable $F$ with respect to $W$ is the random variable 
$$D^W F=\sum_{j=1}^n({\partial_j}f)(\int_0^T h_{1}(s)dW_s,\dots,\int_0^T h_{n}(s)dW_s)h_{j}.$$
 The operator $D^W$ is a linear operator from $L^{2}(\Omega_W)$ into $L^{2}(\Omega_W\times[0,T]),$ closable and unbounded. We will always consider the close extension of $D^W$ and its domain will be
 denoted by $\mathbb{D}^W_{1,2}$.

%

The Skorohod integral with respect to $W$, denoted by $\delta^{W},$ is the adjoint of the derivative operator $D^{W}:\mathbb{D}^W_{1,2} \subset L^{2}\left(\Omega_W\right)\rightarrow L^{2}\left(\Omega_W\times[0,T]\right)$). That is, $u$ is in $Dom \ \delta^W$ if and only if $u\in L^{2}\left(\Omega_W\times[0,T]\right)$) and there exists a random variable $\delta^W(u)\in L^2(\Omega_W)$ satisfying the duality relation 
\begin{equation}\label{eq:rski}
\mathbb{E}_W\left[\int_0^T u_tD_t^WFdt\right]=\mathbb{E}_W\left[\delta^W(u)F\right]\quad\hbox{\rm for every}\quad F\in\mathbb{D}^W_{1,2},
\end{equation}
where $\mathbb{E}_W$ is the expectation with respect to the probability measure $P_W.$

We can extend the last definitions to Hilbert space valued random variables: Let $\mathcal{S}^W(L^2(\Omega_N))$ be the set of all smooth $L^2(\Omega_N)$-random variables of the form 
\begin{equation}\label{eq:smo}
F=\sum_{i=1}^n f_i(\int_0^T h_{1,i}(s)dW_s,\dots,\int_0^T h_{n_i,i}(s)dW_s)G_i,
\end{equation}
where $n\geq 1$, $h_{j,i}\in L^2([0,T])$,  $G_i\in L^2(\Omega_N)$ and $f_i\in \mathcal{C}_{b}^{\infty }(\mathbb{R}^{n_i}),$  for $i\in\{0,\ldots,n\}$ and $j\in\{0\ldots,n_i\}$. The derivative of the random variable $F$ with respect to $W$ is the $L^2(\Omega_N\times[0,T])$-valued random variable
$$D^WF=\sum_{i=1}^n \sum_{j=1}^{n_i}({\partial_j}f_i)(\int_0^T h_{1,i}(s)dW_s,\dots,\int_0^T h_{n_i,i}(s)dW_s) h_{j,i} G_i.$$
 The operator $D^W$ is a linear operator from $L^{2}(\Omega)$ into $L^{2}\left(\Omega\times[0,T]\right),$ closable and unbounded. Moreover it can be iterated defining $D^{W,k}_{t_1,\dots, t_k}F:=D^W_{t_k}\cdots D^W_{t_1}F.$

For any $k,p \geq 1$, we introduce the spaces $\mathbb{D}^W_{k,p}(L^2(\Omega_N))$ as the closure of $\mathcal{S}^{W}(L^2(\Omega_N))$ with respect to the norm 
\begin{equation*}
\left\| F\right\|^p_{W,k,p}:=||\,\, |F|_{L^2(\Omega_N)}\,\, ||_{
L^p(\Omega_W)}^p+\sum_{j=1}^k\left\| \left(\int_{[0,T]^j} |D^{W,j}_z F|_{L^2(\Omega_N)}^2 dz\right)^{\frac{1}{2}}\right\|^p_{L^p(\Omega_W)}.
\end{equation*}

Now the  Skorohod integral with respect to $W$, denoted by $\delta^{W},$ is the adjoint of the derivative operator $D^{W}:\mathbb{D}^W_{1,2}(L^2(\Omega_N)) \subset L^{2}\left(\Omega\right)\rightarrow L^{2}\left(\Omega\times[0,T]\right)$.
That is, $u$ is in $Dom \ \delta^W$ if and only if $u\in L^{2}\left(\Omega\times[0,T]\right)$ and there exists a  random variable $\delta^W(u)\in L^2(\Omega)$ satisfying the duality relation 

\begin{equation}\label{eq:ski}
\mathbb{E}\left[\int_0^T u_tD_t^WFdt\right]=\mathbb{E}\left[\delta^W(u)F\right] \quad\hbox{\rm for every}\quad F\in\mathbb{D}^W_{1,2}(L^2(\Omega_N).
\end{equation}
 The operator 
  $\delta^W$ is an extension of the It\^{o} integral in the sense that the set $L_{a}^{2}(\Omega_W \times [0,T])$ of all square-integrable and adapted processes with respect  to the filtration generated by $X$ is included in $Dom\ \delta^{W}$ and the operator $\delta^{W}$ restricted to $L_{a}^{2}(
  \Omega_W \times [0,T] )$ coincides with the It\^{o} stochastic integral with respect to $W$. For $u\in Dom \ \delta^{W}$ we 
will make use of the notation $\delta^{W}(u)=\int_{0}^{T}u_{t}\delta W_{t}$ and for $u{1\!\!1}_{[0,t]}$  in $Dom\ \delta^W$ we will write $\delta^{W}(u{1\!\!1}_{[0,t]})=\int_{0}^{t}u_{s}\delta W_{s}.$
Note that in (\ref{eq:rski}) and (\ref{eq:ski}) we are using 
$\delta^W$ for the Skorohod integrals defined on $L^2(
\Omega_W\times[0,T])$ and on $L^2(\Omega\times[0,T])$, respectively.
 We
hope that the space will be clear when we use this operator.
 
The following result will be important in  next section.

\begin{lemma}\label{lem:tdi}
Let $F\in\mathbb{D}^W_{1,2}(L^2(\Omega_N))$ and $u\in Dom\ \delta^W\cap L^{2}\left(\Omega\times[0,T]\right)$. Then, for almost
all $\omega''\in\Omega_N$, $F(\cdot,\omega'')\in  \mathbb{D}^W_{1,2}$,  $u(\cdot,\omega'')\in Dom\ \delta^W\cap L^{2}\left(\Omega_W\times[0,T]\right)$,
$$D^WF(\cdot,\omega'')=(D^WF)(\cdot,\omega'')$$
and 
$$\delta^W(u(\cdot,\omega''))=\delta^W(u)(\cdot,\omega'').$$
\end{lemma} 

\noindent\textbf{Remark} Note that left-hand sides of last two
equalities are given by (\ref{eq:rsmo}) and (\ref{eq:rski}), while
right-hand sides are defined via (\ref{eq:smo}) and (\ref{eq:ski}),
 respectively.

\noindent\textbf{Proof of Lemma \ref{lem:tdi}:}
 Let $F\in\mathbb{D}^W_{1,2}(L^2(\Omega_N))$. Then, there is a sequence $\{F_n\in\mathcal{S}^W(L^2(\Omega_N)):n\in
\mathbb{N}\}$ of the form (\ref{eq:smo}) such that $\left\| F_n-F\right\|_{W,1,2}\rightarrow 0$. Hence, the  definition of the canonical L\'evy space, in particular the definition of the probability measure $P$, implies that there is a subsequence $\{n_k: k\in\mathbb{N}\}$ such that, for a.a. $\omega''\in \Omega_N$,
$$||F_{n_k}(\cdot,\omega'')-F(\cdot,\omega'')||_{L^2(\Omega_W)}^2+\left\| \left(\int_{[0,T]} |(D^{W}_z F_{n_k}
(\cdot,\omega''))-(D^{W}_z F)(\cdot,\omega'')|^2 dz\right)^{\frac{1}{2}}\right\|^2_{L^2(\Omega_W)}
\rightarrow 0,$$ 
which gives that the first part of the result is true because $\{F_{n_k}(\cdot,\omega''):k\in \mathbb{N}\}$ is a sequence of the form (\ref{eq:rsmo}).

Finally, let $H\in\mathcal{S}^W$ and $G\in L^2(\Omega_N)$. Then, the duality relation (\ref{eq:ski}) yields
$$\mathbb{E}\left[G\int_0^T u_tD_t^WHdt\right]=\mathbb{E}\left[G\delta^W(u)H)\right].$$
Consequently, using the definition of the probability measure $P$,
for a.a. $\omega''\in\Omega_N$,
$$\mathbb{E}_W\left[\int_0^T u_t(\cdot,\omega'')D_t^WHdt\right]=\mathbb{E}_W\left[\delta^W(u)(\cdot,\omega'')
H(\cdot,\omega'')\right].$$
Thus, from the duality relation (\ref{eq:rski}), the  proof is complete.
\hfill$\square$

\subsection{Anticipative Girsanov Transformations}
 Here, for the convenience of the reader, we recall some basic facts on anticipative Girsanov transformations. 
By Lemma \ref{lem:tdi}, some of these results will be a consequence of the properties of transformations on Wiener space. For a more detailed account on this subject we refer to \cite{Bu,Bu2,bu3}.
Remember that, by the definition of the Canonical L\'evy space, we have that for any $\omega\in\Omega$ there
are $\omega'\in\Omega_W$ and $\omega''\in\Omega_N$ such that $\omega=(\omega',\omega'')$ and viceversa. For $\omega'\in\Omega_W$ and $\omega''\in\Omega_N$, we use the convention $\omega=(\omega',\omega'').$

Given a process $a\in L^2 (\Omega\times [0,T])$, we define the transformation $T_a:\Omega\rightarrow\Omega_W$ as the application  defined by 
$$T_a(\omega',\omega''):=\omega'+\int_0^{\cdot} a_s (\omega',\omega'')ds.$$
Observe that for $\omega''$ fixed, we obtain a transformation on the Wiener space. We say this transformation is absolutely continuous if the measure ${P_W}\circ (T_a(\cdot,\omega''))^{-1}$ is absolutely continuous with respect to $ P_W$, for almost all $\omega''\in\Omega_N$.
Henceforth, we introduce the Cameron-Martin space $CM$, that is, the subspace of absolutely continuous functions of 
$\Omega_W,$ with square-integrable derivatives, endowed with the norm
$$|\omega'|_{CM}:=\left(\int_0^T {\dot\omega}'(t)^2 dt\right)^{\frac{1}{2}}.$$

The following two results are an immediate consequence of \cite{Bu,Bu2,bu3} and Lemma \ref{lem:tdi}:

\begin{proposition}\label{firstinequality}
Let $T^1$ and $T^2$ be  two absolutely continuous transformations associated with processes $a_1$ and $a_2$, respectively, 
$F\in {\mathbb D}^W_{1,2}(L^2(\Omega_N))$ and $\sigma\in L^2([0,T],{\mathbb D}^W_{1,2}(L^2(\Omega_N))).$ Then, for almost all 
$\omega''\in\Omega_N$, we have
\begin{displaymath}|F(T_{a_1}(\omega',\omega''),\omega'')-F(T_{a_2}(\omega',\omega''),\omega'')|\leq || |D^WF|_2 ||_{\infty} |T_{a_1}(\omega',\omega'')-T_{a_2}(\omega',\omega'')|_{CM}\end{displaymath}
and
\begin{eqnarray*}\lefteqn{\left(\int_0^T |\sigma_s (T_{a_1}(\omega',\omega''),\omega'')-\sigma_s (T_{a_2}(\omega',\omega''),\omega'')|^2 ds\right)^{\frac{1}{2}}}\\
&&\leq \left\|\left(\int_0^T \int_0^T |D_r^W \sigma_s|^2 ds dr\right)^{\frac{1}{2}}\right\|_{\infty}|T_{a_1}(\omega',\omega'')-T_{a_2}(\omega',\omega'')|_{CM}. 
\end{eqnarray*}
\end{proposition}
 
\begin{proposition}\label{propbeta}
Let $T_a$ be an absolutely continuous  transformation. Assume 
$a\in L^2 ([0,T];$ ${\mathbb D}^W_{1,2}(L^2(\Omega_N))),$
and let $\sigma\in L^2 ([0,T];{\mathbb D}^W_{1,2}(L^2(\Omega_N)))$ be with 
$\sigma(T_a,\cdot)\in L^2 ([0,T], L^2(\Omega))$ and 
\begin{displaymath}\left\|\left(\int_0^T \int_0^T |D_r^W \sigma_s|^2 ds dr\right)^{\frac{1}{2}}\right\|_{\infty}<\infty.\end{displaymath}
Then, for almost all $\omega''\in\Omega_N$, we get
$\sigma(T_a(\cdot,\omega''),\omega'')\in L^2 ([0,T]
;{\mathbb D}^W_{1,2}),$
\begin{eqnarray*}
D_t^W(\sigma_s(T_a(\omega',\omega''),\omega''))&=&(D_t^W\sigma_s)(T_a(\omega',\omega''),\omega'')\\
&&+\int_0^T (D_r^W\sigma_s)(T_a(\omega',\omega''),\omega'')(D_t^Wa_r)(\omega',\omega'')dr,
\end{eqnarray*}
and
 \begin{eqnarray*}
\int_0^T \sigma_s (T_a(\omega',\omega''),\omega'')\delta W_s&=&\left(\int_0^T \sigma_s \delta W_s\right)(T_a(\omega',\omega''),\omega'')\\
&&-\int_0^T \sigma_s(T_a(\omega',\omega''),\omega'')a_s(\omega',\omega'') ds\\
&&-\int_0^T\int_0^T (D^W_r\sigma_s)(T_a(\omega',\omega''),\omega'')(D^W_sa_r )(\omega',\omega'') drds,
\end{eqnarray*}
for almost all $\omega'\in\Omega_W$.
\end{proposition}

In the remaining of this paper ${\mathbb D}^W_{1,\infty}(L^2(\Omega_N))$ represents the elements $F$ in ${\mathbb D}^W_{1,2}(L^2(\Omega_N))$ such that 
$$||F||_{1,\infty}:=||F||_{\infty}+||\, |D^WF|_2\,||_{\infty}<\infty.$$
Similarly ${\mathbb D}^W_{2,\infty}(L^2(\Omega_N))$ is the family of all the elements in ${\mathbb D}^W_{2,2}(L^2(\Omega_N))
\cap {\mathbb D}^W_{1,\infty}(L^2(\Omega_N))$ such that $D^{W,2}F\in L^{\infty}(\Omega;L^2([0,T]^2)$.

Now, for  $a\in L^2 ([0,T];{\mathbb D}^W_{1,\infty}(L^2(\Omega_N)))$ fixed, we consider two a families of transformations $\{T_t:\Omega\rightarrow \Omega_W: 0\le t \le T\}$  and $\{A_{s,t}:\Omega\rightarrow\Omega_W: 0\le s\le t \le T\}$, which are
the solutions of the equations 
\begin{equation}\label{eq3}
(T_t\, \omega)_{\cdot}=\omega_{\cdot}'+\int_0^{t\wedge\cdot} a_s(T_s\omega,\omega'')\ ds.
\end{equation}
and
\begin{equation}\label{eq3'a}
(A_{s,t}\, \omega)_{\cdot}=\omega'_{\cdot}-\int_{s\wedge\cdot}^{t\wedge\cdot} a_r(A_{r,t} \omega,\omega'')\ dr,
\end{equation}
respectively.

Observe that, for  simplicity of the notation, we not make explicit the dependence on $a$ in these equations. Some of the  properties of the solutions to (\ref{eq3}) and (\ref{eq3'a}) that we need are established in the following result. See \cite{Bu,Bu2,bu3} for its proof.

\begin{proposition}\label{propgamma}
Let $a\in L^2 ([0,T]; {\mathbb D}^W_{1,\infty}(L^2(\Omega_N)))$. Then, there exist two unique families of absolutely continuous transformations $\{T_t, 0\le t\le T\}$ and $\{A_{s,t}: 0\le s\le t \le T\}$ that satisfy  equations (\ref{eq3}) and (\ref{eq3'a}), respectively. Moreover, for each $s,t\in[0,T]$, $s<t$, $A_{s,t}=T_sA_t$, with $A_t=A_{0,t}$, $T_t$ is invertible with  inverse $A_t$ and  $a_{\cdot}(T_{\cdot}(\ast,\omega''),\omega'')\in L^2([0,T]; \mathbb{D}^W_{1,\infty})$, for a.a. $\omega''\in\Omega_N$.
\end{proposition}

In relation to the transformation $A_{s,t},$ we have the following lemma that will be useful for our purposes. 

\begin{lemma}\label{uniformcm}
Let $a\in L^2 ([0,T];{\mathbb D}^W_{1,\infty}(L^2(\Omega_N))).$ Then, for any $u\leq s \leq t$,  we have 
$$|A_{u,t}\omega-A_{u,s}\omega|^2_{CM}\leq 2\left(\int_s^t ||a_r||^2_{\infty} dr\right)\exp\left\{2\int_0^T ||\,
 |D^Wa_r|_2^2\, ||_{\infty} dr\right\}.$$
\end{lemma}

\noindent\textbf{Remark} Note $A_{u,t}$ is continuous in $t$ with respect the CM-norm, uniformly on $u$. 

\medskip

\noindent{\bf Proof of Lemma \ref{uniformcm}}
Let $u\leq s\leq t.$ Then, by Propositions \ref{firstinequality} and \ref{propgamma} we have 
\begin{eqnarray*} 
|A_{u,s}\omega-A_{u,t}\omega|^2_{CM}&=&\left|\int_{u\wedge\cdot}^{s\wedge\cdot} a_r (A_{r,s}\omega,\omega'')dr-\int_{u\wedge\cdot}^{t\wedge\cdot} a_r (A_{r,t}\omega,\omega'')dr\right|_{CM}^2\\
&=&\int_0^T |{1\!\!1}_{(u,s]}(r) a_r (A_{r,s}\omega,\omega'')-{1\!\!1}_{(u,t]} (r)a_r(A_{r,t}\omega,\omega'')|^2 dr\\
&\leq&2\int_s^t |a_r(A_{r,t}\omega,\omega'')|^2dr+2\int_u^s |a_r(A_{r,s}\omega,\omega'')-a_r(A_{r,t}\omega,\omega'')|^2 dr\\
&\leq&2\int_s^t ||a_r||^2_{\infty} dr+2\int_u^s ||\,
 |D^Wa_r|_2^2\, ||_{\infty} |A_{r,s}\omega-A_{r,t}\omega|^2_{CM} dr.
\end{eqnarray*}
So, using Gronwall's lemma, we obtain
$$|A_{u,t}\omega-A_{u,s}\omega|^2_{CM}\leq 2\left(\int_s^t ||a_r||^2_{\infty} dr\right)\exp\left\{2\int_u^s ||\, 
|D^Wa_r|_2^2\, ||_{\infty} dr\right\}$$
which implies the result holds. 
\hfill{$\square$}

\medskip

\noindent\textbf{Remark} In \cite{Bu,Bu2,bu3}, Buckdahn has proven 
that both inequalities in Proposition \ref{firstinequality} hold
only for almost all $\omega'\in\Omega_W$. But, by Fubini theorem,
 it
 is not difficult to see that, in this case, the inequality in
 Lemma \ref{uniformcm} is satisfied for a.a. $\omega\in\Omega$.

\medskip 

To finish this subsection we give some results related to the densities  of the transformations
$\{T_t:\Omega\rightarrow \Omega_W: 0\le t \le T\}$  and $\{A_{s,t}:\Omega\rightarrow\Omega_W: 0\le s\le t \le T\}$. 
Now, let $F\in L^{\infty}(\Omega)$ and $a$ as in Proposition \ref{propgamma}. One of our main tools in the proof of the existence and uniqueness of the solution to equation (\ref{eintr2}) are the equalities, proven by Buckdahn \cite{Bu,Bu2,bu3},
\begin{equation}
\mathbb{E}\left[F(A_{s,t}\omega,\omega'') L_{s,t}(\omega)\right]= \mathbb{E}\left[F\right]\label{egirsanov1}
\end{equation}
and
\begin{equation}
\mathbb{E}\left[F(A_{s,t}\omega,\omega'') \right]= \mathbb{E}\left[F \mathcal{L}_{s,t}\right],\label{egirsanov2}
\end{equation}
where  
\begin{eqnarray}\label{edensitya}
L_{s,t}(\omega)&=&\exp\left\{\int_{s}^t a_r(A_{r,t}\omega,\omega'')\delta W_r-\frac12 \int_{s}^t a_r^2(A_{r,t}\omega,\omega'')dr\right.\nonumber\\
&&\quad\quad\quad\left.-\int_{s}^t\int_r^t (D^W_ua_r)(A_{r,t}\omega,\omega'')D^W_r[a_u(A_{u,t}\omega,\omega'')]
dudr\right\}     \end{eqnarray}
is the density of $A_{s,t}^{-1}$ and 
\begin{eqnarray}\label{edensityb}
\mathcal{L}_{s,t}(\omega)&=&\exp\left\{-\int_{s}^t a_r(T_tA_r
\omega,\omega'')\delta W_r-\frac12 \int_{s}^t a_r^2(T_tA_r\omega,\omega'')dr\right.\nonumber\\
&&\quad\quad\quad\left.-\int_{s}^t\int_s^r (D_u^Wa_r)(T_tA_r\omega,\omega'')D_r^W[a_u(T_tA_u\omega,\omega'')]dudr\right\}.
\end{eqnarray} 
Finally, we have that, in this case,
\begin{eqnarray}
L_{s,t}(\omega)&=& \mathcal{L}^{-1}_{s,t}(A_{s,t}\omega,\omega''),\label{erelacion1}\\
L_{0,t}(\omega)&=&L_{0,s}(A_{s,t}\omega,\omega'')
 L_{s,t}(\omega),\quad 0\le s\le t\le T.\label{relacion2}
\end{eqnarray}
These two relations  can be proved as consequence of the equalities (\ref{egirsanov1}) and (\ref{egirsanov2}). Indeed, 
\begin{displaymath}
\mathbb{E}\left[F(A_{s,t}\omega,\omega'') L_{s,t}(\omega)\right]= \mathbb{E}\left[F\right]
=\mathbb{E}\left[F \mathcal{L}_{s,t} \mathcal{L}^{-1}_{s,t}\right]=\mathbb{E}\left[F(A_{s,t}\omega,\omega'')
\mathcal{L}^{-1}_{s,t}(A_{s,t}\omega,\omega'')\right],
\end{displaymath}
and 
\begin{displaymath}
\mathbb{E}\left[F(A_{t}\omega,\omega'') L_{0,t}(\omega)\right]= \mathbb{E}\left[F(A_{s}\omega,\omega'') L_{0,s}(\omega)\right]
=\mathbb{E}\left[F(A_{t}\omega,\omega'') L_{0,s}(A_{s,t}\omega,\omega'') L_{s,t}(\omega)\right].
\end{displaymath}

\subsection{The anticipative linear stochastic differential equation on canonical Wiener space}

On the canonical Wiener space, Buckdahn \cite{Bu,Bu2,bu3} has studied equation (\ref{eintr2}) via the anticipating Girsanov transformations (\ref{eq3}) and (\ref{eq3'a}).  Namely, he considers the linear stochastic differential equation
\begin{equation}\label{eintr7}
Z_t=Z_{0}+ \int_{0}^t h_s Z_s\ ds +  \int_{0}^t a_s Z_s\ \delta W_s,\qquad t\in [0,T],
\end{equation}
and state the following result:
\begin{theorem}\label{pintro1}
Assume $a\in L^2([0,T], \mathbb{D}^W_{1,\infty})$, $h\in L^1([0,T], L^\infty(\Omega))$ and $Z_{0}\in L^\infty(\Omega)$. 
Then, the process $Z=\{Z_t:t\in[0,T]\}$ defined by
\begin{equation}\label{eintr8}
Z_t:=Z_{0}(A_{0,t})\exp\left\{\int_{0}^t h_s (A_{s,t})\
ds\right\}\ L_{0,t}
\end{equation}
belongs to $L^1(\Omega\times [0,T])$ and is a global solution of (\ref{eintr7}). Conversely, if $Y\in L^1(\Omega\times [0,T])$ is a global solution  of (\ref{eintr7}) and, if, moreover,
 $a, h \in L^\infty(\Omega\times [0,T])$ and 
 $D^Wa \in L^\infty(\Omega\times [0,T]^2)$, then $Y$ is of the form (\ref{eintr8})
for a.e. $0\le t\le T$.
\end{theorem}

Moreover we need the following proposition on the continuity of $Z$, whose proof is given in the Appendix (see Section
 \ref{sec:ape})
because it is too long and technical.
\begin{theorem}\label{continuity}
Assume $Z_{0}\in {\mathbb D}^W_{1,\infty}$, $h\in L^1([0,T], \mathbb{D}^W_{1,\infty})$ and that, for some $p>2$, 
$$a\in L^{2p}([0,T], \mathbb{D}^W_{1,\infty})\cap L^{2}([0,T], \mathbb{D}^W_{2,\infty}).$$ 
Then, $Z$ given by (\ref{eintr8}) has continuous trajectories a.s. 
\end{theorem}

\section{Two processes with jumps}
\setcounter{equation}{0}
\setcounter{theorem}{0}

 In the sequel we use the following hypothesis on the coefficients:

\begin{itemize} 
\item[\textbf{(H1)}] Assume that $X_0\in {\mathbb D}^W_{1,\infty}(L^\infty(\Omega_N))$, 
$ b, v_{\cdot-}(y)\in L^1([0,T], {\mathbb D}^W_{1,\infty}(L^\infty(\Omega_N))$, for all 
$y\in {\R}_0$. Moreover, there exists $p>2$ such that 
$$a\in L^2([0,T],\mathbb{D}^W_{2,\infty}(L^\infty(\Omega_N))\cap L^{2p}([0,T], {\mathbb D}^W_{1,\infty}(L^\infty(\Omega_N)).$$
\item[\textbf{(H2)}] There exist a positive function $g\in L^2(\mathbb{R}_0,\nu) \cap
L^1(\mathbb{R}_0,\nu)$ such that 
\begin{displaymath}
\left|v_{s-}(y, \omega)\right| \le g(y),\qquad \textrm{uniformly on}\ \omega\ \textrm{and} \ s,
\end{displaymath}
and
\begin{displaymath}
\lim_{|y|\rightarrow 0} g(y)=0.
\end{displaymath}
\item[\textbf{(H3)}]
The function $g$ satisfies $\int_{{\mathbb R}_0} (e^{g(y)}-1)\nu(dy)<\infty.$
\item[\textbf{(H4)}]
The function $g$ satisfies $\int_{{\mathbb R}_0} (e^{2g(y)}-1)\nu(dy)<\infty.$
\end{itemize}

\noindent {\bf Remark:} As an example, observe that the following function is in $L^1({\mathbb R}_0,\nu)\cap L^2({\mathbb R}_0,\nu)$
and is such that (H3) and (H4) hold, and $\lim_{|y|\rightarrow 0}
g(y)=0$.
\begin{displaymath}
g(y)=\left\{\begin{array}{ll}
 k_1(\beta) y^2, & y\in (-\beta,\beta),\\
 k_2(\beta), & y\in (-\beta,\beta)^c,
\end{array}\right.
\end{displaymath}
where $\beta\in (0,1)$ and $k_1(\beta)$ and $k_2(\beta)$ are positive constants.

\bigskip

Given $\varepsilon>0,$ set
\begin{eqnarray}
X_t^\varepsilon&=&X_0(A_{0,t})\exp\left\{\int_{0}^t b_s (A_{s,t})\ ds\right\}\  L_{0,t}\ \prod_{s\le t, \varepsilon<|y|}
\Big[1+v_{s-}(y,A_{s,t}) \Delta N(s,y)\Big]\nonumber\\
&&\times \exp\left\{-\int_0^t \int_{\{|y|>\varepsilon\}} v_{s-}(y, A_{s,t}) \ \nu(dy) ds\right\}\label{eintr9}.
\end{eqnarray}
Notice that this process can also be written as follows
\begin{displaymath}
X_t^\varepsilon=X_0( A_{0,t})\exp\left\{\int_{0}^t b_s^\varepsilon
(A_{s,t})\ ds\right\}\  L_{0,t}\ \prod_{i=1}^{N_t^\varepsilon}
\Big[1+v_{\tau_i^\varepsilon-}(y_i^\varepsilon,
A_{\tau_i^\varepsilon,t})\Big],\end{displaymath} 
where $b_s^\varepsilon (\omega):=b_s(\omega)-\int_{\{|y|>\varepsilon\}} v_{s-}(y,\omega)\nu(dy),$ $\{\tau_i^{\varepsilon}, i\ge 1\}$ are the jumps times whose jumps size are
greater than $\varepsilon$, $y_i^\varepsilon$ denotes the amplitude of jump $\tau_i^\varepsilon$, and $N_t^\varepsilon$ is the number of jumps before $t$, with size bigger that $\varepsilon$.

\medskip

\begin{proposition}\label{almostsurely}
Assume (H1) and (H2) hold. For each $t\in [0,T]$, the process $X_t^\varepsilon$, defined in (\ref{eintr9}), converges almost surely to
\begin{eqnarray}
X_t&=&X_0(A_{0,t})\exp\left\{\int_{0}^t b_s(A_{s,t})\ ds\right\}\  L_{0,t}\ \exp\left\{-\int_0^t \int_{\mathbb{R}_0} v_{s-}(y,
A_{s,t}) \ \nu(dy)\ ds\right\}\nonumber\\
&&\times\prod_{s\le t, y\in \mathbb{R}_0}\Big[1+v_{s-}(y,A_{s,t}) \Delta N(s,y)\Big]\label{eintr9bis}.
\end{eqnarray} 
\end{proposition}

\noindent {\bf Remark:} In the proof of this result we will see
that the representation
\begin{eqnarray*}
X_t&=&X_0(A_{0,t})\exp\left\{\int_{0}^t b_s(A_{s,t})\ ds\right\}\  L_{0,t}\ \exp\left\{\int_0^t \int_{\mathbb{R}_0} v_{s-}(y,
A_{s,t}) \ d\tilde N(s,y)\right\}\nonumber\\
&&\times\prod_{s\le t, y\in \mathbb{R}_0}\Big[1+v_{s-}(y,A_{s,t}) \Delta N(s,y)\Big] e^{-v_{s-}(y,A_{s,t})\Delta N(s,y)}
\end{eqnarray*} 
also holds. We observe that the stochastic integral with respect to
${\tilde N}$ is pathwise defined.

\medskip

\noindent {\bf Proof of Proposition \ref{almostsurely}:} First of all, the hypotheses on $X_0$ and $b$ yield
\begin{displaymath}
\left|X_0(A_{0,t})\exp\left\{\int_{0}^t b_s
(A_{s,t})\ ds\right\}\right|\le C.
\end{displaymath}
Secondly, as the factor $L_{0,t}$ is a density, it is finite a.s. So it  remains to see the convergence of the following
quantities:
\begin{eqnarray*}
M_1&=&\exp\left\{\int_0^t \int_{|y|>\varepsilon} v_{s-}(y,
A_{s,t}) \ d\tilde N(s,y)\right\},\\
M_2&=&\prod_{s\le t,  \varepsilon<|y|}
\Big[1+v_{s-}(y,
A_{s,t}) \Delta N(s,y)\Big] e^{-v_{s-}(y,
A_{s,t})\Delta N(s,y)}.
\end{eqnarray*}
Using the relation $d\tilde N(t,y)=d N(t,y)- \nu(dy) dt$ and (H2), we have
\begin{equation}\label{emexico2}\begin{array}{l}
\displaystyle
\int_0^t \int_{\mathbb{R}_0} \left|v_{s-}(y,
A_{s,t})\right| \ d N(s,y)+ \int_0^t \int_{\mathbb{R}_0} \left|v_{s-}(y,
A_{s,t})\right| \ \nu(dy) ds\\[3mm]
 \displaystyle\qquad \le 
\int_0^t \int_{\mathbb{R}_0} g(y) \ d N(s,y)+ \int_0^t \int_{\mathbb{R}_0} g(y) \ \nu(dy) ds\\[3mm]
 \displaystyle\qquad =
\int_0^t \int_{\mathbb{R}_0} g(y) \ d\tilde N(s,y)+ 2 \int_0^t \int_{\mathbb{R}_0} g(y) \ \nu(dy) ds.\end{array}
\end{equation}
This quantity is finite a.s. because (H2) implies
\begin{displaymath}
\mathbb{E}\left[\left(\int_0^t \int_{\mathbb{R}_0} g(y) \ d\tilde N(s,y)\right)^2\right]
=\int_0^t \int_{\mathbb{R}_0} g(y)^2 \ \nu(dy) ds<\infty.
\end{displaymath}
 Then, $M_1$ converges a.s., as $\varepsilon\rightarrow 0$, to
 $\exp\{\int_0^t\int_{{\mathbb R}_0}v_{s-}
 (y,A_{s,t})d{\tilde N}(s,y)\}$.

\medskip

On the other hand, for any constant  $c>0$,
\begin{displaymath}
M_2=M_{2,1}\times  M_{2,2},
\end{displaymath}
with
\begin{eqnarray*}
M_{2,1}&=&\prod_{s\le t,  \varepsilon<|y|<c\vee\varepsilon}
\Big[1+v_{s-}(y,
A_{s,t}) \Delta N(s,y)\Big] e^{-v_{s-}(y,
A_{s,t})\Delta N(s,y)},\\
M_{2,2}&=&\prod_{s\le t, c\vee \varepsilon\le|y|}
\Big[1+v_{s-}(y,
A_{s,t}) \Delta N(s,y)\Big] e^{-v_{s-}(y,
A_{s,t})\Delta N(s,y)}.
\end{eqnarray*}
$M_{2,2}$ is well-defined and converges as $\varepsilon\rightarrow 
0$ to
$$\prod_{s\le t, c\le|y|}
\Big[1+v_{s-}(y,
A_{s,t}) \Delta N(s,y)\Big] e^{-v_{s-}(y,
A_{s,t})\Delta N(s,y)},$$
because it is a product of a finit number of factors.
To deal with $M_{2,1}$ we use the following argument:
 Hypothesis (H2) implies that for small enough $y$,  $\left|v_{s-}(y,\omega)\right| \le\frac12$.
Then, choosing $c>0$ such that
$|g(y)|<\frac12$, for $|y|<c$,
\begin{displaymath}
\log M_{2,1}= \sum_{s\le t,\varepsilon<|y|<c}\left[\log\left(1+v_{s-}(y,
A_{s,t}) \Delta N(s,y)\right)-v_{s-}(y,
A_{s,t})\Delta N(s,y)\right],
\end{displaymath}
and this series is absolutely convergent since
\begin{displaymath}\begin{array}{l}
\displaystyle \sum_{s\le t,\varepsilon<|y|<c}\left|\log\left(1+v_{s-}(y,
A_{s,t}) \Delta N(s,y)\right)-v_{s-}(y,
A_{s,t})\Delta N(s,y)\right|\\[2mm]
\displaystyle \qquad \qquad \le \sum_{s\le t,0<|y|<c} \left[v_{s-}(y,A_{s,t}) \Delta N(s,y)\right]^2 \le \sum_{s\le t,0<|y|<c} \frac12 \left|v_{s-}(y,A_{s,t}) \Delta N(s,y)\right|\\[4mm]
\displaystyle \qquad \qquad \le  \frac12 \sum_{s\le t,0<|y|<c} g(y) \Delta N(s,y).\end{array}
\end{displaymath} So, $M_{2,1}$ also converges
as $\varepsilon\rightarrow 0$ 
since  $\int_0^t \int_{\mathbb{R}_0} g(y) \ d N(s,y)$ is finite by (\ref{emexico2}).

We can conclude that the processes $X_t^\varepsilon$ and $X_t$ are well-defined and $X_t^\varepsilon$ converges a.s. to $X_t$ when $\varepsilon$ goes to zero.
\hfill{$\square$}

\medskip

\begin{proposition}\label{ela1}
Assume (H1), (H2) and (H3). Then, the processes $X^\varepsilon$ and $X$ belong to $L^1(\Omega\times [0,T])$ and
\begin{equation}
X_t^\varepsilon\tresp\hbox to 12mm{\rightarrowfill}_{\varepsilon\rightarrow 0}^{} X_t,
\label{emexico3}
\end{equation}
in $L^1(\Omega\times [0,T]).$
\end{proposition}

\noindent {\bf Proof:} Since
 $X_0$ and $b$ are bounded and (H2) is true,  it is immediate to check that
\begin{eqnarray*}
\left|X_t^\varepsilon\right|&=&
\left|X_0(A_{0,t})\right| \exp\left\{\int_{0}^t b_s
(A_{s,t})\ ds\right\}\  L_{0,t}\ \left|\prod_{s\le t, \varepsilon<|y|}
\Big[1+v_{s-}(y,
A_{s,t}) \Delta N(s,y)\Big]\right|\\
&&\times \exp\left\{-\int_0^t \int_{\{|y|>\varepsilon\}} v_{s-}(y,
A_{s,t}) \ \nu(dy) ds\right\}\\
&\le & C L_{0,t}\ \prod_{s\le t, \varepsilon<|y|}
\Big[1+g(y) \Delta N(s,y)\Big]\ \exp\left\{\int_0^t \int_{\{|y|>\varepsilon\}} g(y) \ \nu(dy) ds\right\}\\
&\le & C L_{0,t} \prod_{s\le t, \varepsilon<|y|}
\Big[1+g(y) \Delta N(s,y)\Big],
\end{eqnarray*}
due to  $g\in L^1( \mathbb{R}_0,\nu)$. Moreover, using that $1+x\le e^x$ for $x>0$, we get
\begin{displaymath}
\left|X_t^\varepsilon\right|\le C L_{0,t}  \exp\left\{\sum_{s\le t, \varepsilon<|y|}  g(y) \Delta N(s,y)\right\}
\le C L_{0,t}  \exp\left\{\int_0^T\int_{{\mathbb R}_0} g(y)dN(s,y)
\right\}.
\end{displaymath}
Finally, the result follows from Proposition \ref{almostsurely}, the 
dominated convergence theorem and from the fact that $L_{0,t}
\exp\left\{\int_0^T\int_{{\mathbb R}_0} g(y)dN(s,y)
\right\}\in L^1(\Omega)$, which is a consequence of Cont and Tankov
\cite{CT04} (Proposition 3.6). Indeed, remember that
$E_W(L_{0,t})=1$.
\hfill{$\square$}

\section{Existence and uniqueness of solution of the approximated
 equation} 
\setcounter{equation}{0}
\setcounter{theorem}{0}

The goal of this section is to prove the following theorem.

\begin{theorem}\label{terss}
Assume (H1), (H2) and (H3). Also assume  that $a$, $b$ and $v_{\cdot -}
(y),$ for any $y\in {\mathbb R}_0,$ belong to $L^\infty(\Omega\times [0,T])$ and 
$D^Wa$ belongs to $L^\infty(\Omega\times [0,T]^2).$ Then, the process $X^\varepsilon$ defined in 
(\ref{eintr9}) is the unique solution 
in $L^1(\Omega\times [0,T])$
 of 
\begin{equation}\label{enewversion}
X_t^\varepsilon= X_0+
\int_0^t b_s X_s^\varepsilon ds+\int_0^t a_s X_s^\varepsilon \delta W_s+
 \int_0^t \int_{\{|y|>\varepsilon\}} v_{s-}(y) X_{s-}^\varepsilon\ d \tilde N(s,y).
\end{equation}
\end{theorem}

\medskip

\noindent {\bf Remark}. Note that Equation (\ref{eintr3}) can be rewritten as an equation of the form   (\ref{enewversion}). 

\medskip

\noindent {\bf Proof of Theorem \ref{terss}:} The proof is divided into three steps.

\medskip

\noindent {\bf Step 1.}  We first analyze the continuity of
$X^\varepsilon$. As we have seen before
\begin{displaymath}
X_t^\varepsilon=X_0( A_{0,t})\exp\left\{\int_{0}^t b_s^\varepsilon
(A_{s,t})\ ds\right\}\  L_{0,t}\ \prod_{i=1}^{N_t^\varepsilon}
\Big[1+v_{\tau_i^\varepsilon-}(y_i^\varepsilon,
A_{\tau_i^\varepsilon,t})\Big].
\end{displaymath} 
Observe  that under the hypotheses of the theorem, $X_0( A_{0,t})\exp\left\{\int_{0}^t b_s^\varepsilon(A_{s,t})\ ds\right\}\  L_{0,t}$ is continuous on $t$, for a.a. $\omega''\in\Omega_N$,
as a consequence of Theorem \ref{continuity}. On other hand, $\omega''-$a.s., $\prod_{i=1}^{N_t^\varepsilon}[1+v_{\tau_i^\varepsilon-}(y_i^\varepsilon,
A_{\tau_i^\varepsilon,t})]$ is a finite product with all the terms well defined and continuous on $t$ as a consequence of (H1),
Lemma \ref{uniformcm} and Proposition \ref{firstinequality}. 
So, $X^\varepsilon$ has a.s. right continuous trajectories with left limits. Moreover recall that by Proposition \ref{ela1} the process $X^\varepsilon$ 
belongs to $L^1(\Omega\times [0,T]).$ 

\medskip

\noindent {\bf Step 2.} We now prove that $X^\varepsilon$ is a 
solution of (\ref{enewversion}). 
 
Assume that $G$ is an element of the set of 
$L^2(\Omega_N)$-smooth Wiener functionals described as 
 $G=\sum_{i=1}^n H_i Z_i$, with
$H_i\in \mathcal{S}^W$ and $Z_i\in L^2(\Omega_N)$. Denote
\begin{displaymath}
\Phi^\varepsilon_s(\omega)=\prod_{r\le s, \varepsilon<|y|}
\Big[1+v_{r-}(y,
T_r) \Delta N(r,y)\Big].
\end{displaymath}
Due to Girsanov's theorem (\ref{egirsanov1}) and that $A_{r,t}=T_r\circ A_t$, we have
\begin{eqnarray*}
\mathbb{E}\bigg[\int_0^t a_s X_s^\varepsilon  D_s^WG ds\bigg]
&=&\mathbb{E}\bigg[\int_0^t a_s X_0(A_{0,s}) \exp\bigg\{\int_{0}^s b_r^\varepsilon (A_{r,s}) dr\bigg\}  L_{0,s} \Phi^\varepsilon_s(A_{s})  D_s^WG ds\bigg]\\
&=&\mathbb{E}\bigg[\int_0^t a_s(T_{s}) X_0 \exp\bigg\{\int_{0}^s b_r^\varepsilon
(T_{r}) dr\bigg\} \Phi^\varepsilon_s\ (D_s^WG)(T_s) ds\bigg].
\end{eqnarray*}
Lemma 2.2.4 in \cite{Bu} shows that $\frac{d}{ds} G(T_s)=a_s(T_s)
 (D_s^W G)(T_s)$. So,
\begin{displaymath}
\mathbb{E}\bigg[\int_0^t a_s X_s^\varepsilon  D_s^WG ds\bigg]
= \mathbb{E}\bigg[\int_0^t\left(\frac{d}{ds} G(T_s)\right) X_0 \exp\bigg\{\int_{0}^s b_r^\varepsilon
(T_{r}) dr\bigg\} \Phi^\varepsilon_s  ds\bigg].
\end{displaymath}
Since $\int_0^t \int_{\{|y|>\varepsilon\}}  dN(s,y)<\infty$ a.s., we have 
\begin{displaymath}
\mathbb{E}\bigg[\int_0^t a_s X_s^\varepsilon  D_s^WG ds\bigg]
= \sum_{i=1}^\infty \mathbb{E}\bigg[\int_{\tau_{i-1}^\varepsilon \wedge t}^{\tau_{i}^\varepsilon \wedge t} \left( \frac{d}{ds} G(T_s)\right) X_0 \exp\bigg\{\int_{0}^s b_r^\varepsilon (T_{r}) dr\bigg\} 
\Phi^\varepsilon_{\tau_{i-1}^\varepsilon \wedge t} ds \bigg].
\end{displaymath}
Then, integration by parts implies
\begin{eqnarray*}
\mathbb{E}\bigg[\int_0^t a_s X_s^\varepsilon  D_s^WG ds\bigg]&=&\sum_{i=1}^\infty \mathbb{E}\bigg[G(T_{\tau_{i}^\varepsilon \wedge t}) X_0 \exp\bigg\{\int_{0}^{\tau_{i}^\varepsilon \wedge t} b_r^\varepsilon(T_{r}) dr\bigg\} \Phi^\varepsilon_{\tau_{i-1}^\varepsilon \wedge t}\\
&&- G(T_{\tau_{i-1}^\varepsilon \wedge t}) X_0 \exp\bigg\{\int_{0}^{\tau_{i-1}^\varepsilon \wedge t} b_r^\varepsilon (T_{r}) dr\bigg\} \Phi^\varepsilon_{\tau_{i-1}^\varepsilon \wedge t}\\
&&-\int_{\tau_{i-1}^\varepsilon \wedge t}^{\tau_{i}^\varepsilon \wedge t} G(T_s) X_0 b_s^\varepsilon(T_{s})\exp\bigg\{\int_{0}^s b_r^\varepsilon (T_{r}) dr\bigg\}\Phi^\varepsilon_{\tau_{i-1}^\varepsilon \wedge t} ds\bigg].
\end{eqnarray*}
Using that $\Phi^\varepsilon_{\tau_{i}^\varepsilon}
=\Phi^\varepsilon_{\tau_{i-1}^\varepsilon}
(1+v_{\tau_{i}^\varepsilon -}
(y_i^\varepsilon, T_{\tau_{i}^\varepsilon}))$,
we have that  the previous quantity is equal to
\begin{displaymath}
\begin{array}{l}
\displaystyle \sum_{i=1}^\infty \mathbb{E}\bigg[G(T_{\tau_{i}^\varepsilon \wedge t}) X_0\exp\bigg\{\int_{0}^{\tau_{i}^\varepsilon \wedge t} b_r^\varepsilon(T_{r}) dr\bigg\} \Phi^\varepsilon_{\tau_{i}^\varepsilon \wedge t}\\[3mm]
\displaystyle \qquad - G(T_{\tau_{i-1}^\varepsilon \wedge t}) X_0 \exp\bigg\{\int_{0}^{\tau_{i-1}^\varepsilon \wedge t} b_r^\varepsilon
(T_{r}) dr\bigg\} \Phi^\varepsilon_{\tau_{i-1}^\varepsilon \wedge t}\bigg]\\[3mm]
\displaystyle - \sum_{i; \tau_i\le t} \mathbb{E}\bigg[
G(T_{\tau_{i}^\varepsilon \wedge t}) X_0 \exp\bigg\{\int_{0}^{\tau_{i}^\varepsilon \wedge t} b_r^\varepsilon
(T_{r}) dr\bigg\} v_{\tau_{i}^\varepsilon \wedge t-}(y_i^\varepsilon, T_{\tau_{i}^\varepsilon \wedge t})
\Phi^\varepsilon_{\tau_{i-1}^\varepsilon \wedge t}\bigg]\\[3mm]
\displaystyle -\sum_{i=1}^\infty  \mathbb{E}\bigg[ \int_{\tau_{i-1}^\varepsilon \wedge t}^{\tau_{i}^\varepsilon \wedge t} G(T_s) X_0 b_s^\varepsilon(T_{s})\exp\bigg\{\int_{0}^s b_r^\varepsilon(T_{r}) dr\bigg\} 
\Phi^\varepsilon_{\tau_{i-1}^\varepsilon \wedge t} ds \bigg].\end{array}
\end{displaymath}
Taking into account that the two first summands form a telescopic series, Girsanov's theorem (\ref{egirsanov1}) and (\ref{eintr9}) imply that
\begin{eqnarray*}
\mathbb{E}\bigg[\int_0^t a_s X_s^\varepsilon  D_s^WG ds\bigg]&=&\mathbb{E}\bigg[G(T_{t}) X_0 \exp\bigg\{\int_{0}^{t} b_r^\varepsilon
(T_{r}) dr\bigg\} \Phi^\varepsilon_{t}-GX_0\bigg]\\
&&- \sum_{i; \tau_i\le t} \mathbb{E}\bigg[G X_0(A_{\tau_{i}^\varepsilon \wedge t}) \exp\bigg\{\int_{0}^{\tau_{i}^\varepsilon \wedge t} b_r^\varepsilon
(A_{r,\tau_{i}^\varepsilon \wedge t}) dr\bigg\} L_{0,\tau_{i}^\varepsilon \wedge t}\\
&&\times\ v_{\tau_{i}^\varepsilon \wedge t-}(y_i^\varepsilon)\Phi^\varepsilon_{\tau_{i-1}^\varepsilon \wedge t}(A_{\tau_{i}^\varepsilon \wedge t})\bigg]
-\sum_{i=1}^\infty \mathbb{E}\bigg[ \int_{\tau_{i-1}^\varepsilon \wedge t}^{\tau_{i}^\varepsilon \wedge t} G b_s^\varepsilon X_s^\varepsilon ds\bigg]\\
&=&\mathbb{E}\big[G  X_t^\varepsilon -G X_0\big] - \mathbb{E}\bigg[G\int_0^t \int_{\{|y|>\varepsilon\}} v_{s-}(y) X_{s-}^\varepsilon\ d N(s,y)\bigg]\\
&&-\mathbb{E}\bigg[G\int_0^t b_s^\varepsilon X_s^\varepsilon\ ds\bigg].
\end{eqnarray*}
So,
\begin{equation}\label{emexico7}
\mathbb{E}\bigg[\int_0^t a_s X_s^\varepsilon  D_s^WG ds
\bigg]=\mathbb{E}\bigg[G\Big( X_t^\varepsilon - X_0 - \int_0^t \int_{\{|y|>\varepsilon\}} v_{s-}(y) X_{s-}^\varepsilon\ d \tilde N(s,y)-\int_0^t b_s X_s^\varepsilon ds\Big) \bigg].
\end{equation}
That means that $X^\varepsilon$ defined in (\ref{eintr9}) is solution of (\ref{enewversion}).

\medskip

\noindent {\bf Step 3.} Now we prove the uniqueness of the solution to (\ref{enewversion}). We argue it by induction with respect to the jump
times $\tau_i^\varepsilon$. Notice that if
 $t\in [0,\tau_1^\varepsilon)$, by Theorem \ref{pintro1},
  there exists a unique solution. We now suppose that
   $t \in [\tau_1^\varepsilon,\tau_2^\varepsilon)$.
    Assume that $Y^\varepsilon$  is a solution of the 
    stochastic differential equation (\ref{enewversion})
     such that $Y^\varepsilon\in L^1(\Omega\times[0,T])$ and
      $a_{\cdot}
Y^\varepsilon_{\cdot}\ \1_{[\tau_1^\varepsilon,
\tau_2^\varepsilon)}(\cdot)$ is Skorohod integrable.  For any $t\in [\tau_1^\varepsilon,
\tau_2^\varepsilon)$, $Y^\varepsilon_t$ satisfies
\begin{equation}\label{eqr0mex}
Y_t^\varepsilon=X^\varepsilon_{\tau_1^\varepsilon}
+ \int_{\tau_1^\varepsilon}^tb_s^\varepsilon Y^\varepsilon_s ds + \int_{\tau_1^\varepsilon}^t a_s Y^\varepsilon_s \delta W_s,
\end{equation}
where by Step 1, we can write
\begin{equation}\label{eufff}
X_{\tau_1^\varepsilon}^\varepsilon=X_0( A_{0,\tau_1^\varepsilon})
\Big[1+v_{\tau_1^\varepsilon-}(y_1^\varepsilon)\Big]
\exp\left\{\int_{0}^{\tau_1^\varepsilon} b_s^\varepsilon
(A_{s,\tau_1^\varepsilon})\ ds\right\}\  L_{0,\tau_1^\varepsilon}.
\end{equation} 
Note that Lemma \ref{lem:tdi} implies that, for a.a. $\omega''\in
\Omega_N$, $a(\cdot,\omega'')Y(\cdot,\omega'')\ \1_{]\tau_i^\varepsilon,t]}
\in$ Dom$\delta^W$, and $a(\cdot,\omega'')$, $X_0(\cdot,\omega'')$,
 $b(\cdot,\omega'')$, $v(y,\cdot,\omega'')$ satisfy (H1) when we write
 ${\mathbb D}^W_{1,\infty}$ and ${\mathbb D}^W_{2,\infty}$ instead of
 ${\mathbb D}^W_{1,\infty}(L^\infty(\Omega_N))$ and 
 ${\mathbb D}^W_{2,\infty}(L^\infty(\Omega_N))$, respectively. Now we fix
 a such $\omega''$ and in the following calculations we avoid write
 it to simplify the notation.
 
 For $\rho>0$, there is a sequence $\{a^n:n\in{\mathbb N}\}$ of smooth functionals of the form
 $a^n=\sum_{i=1}^{m_n}F_{i,n}h_{i,n}$, with $F_{i,n}\in{\cal S}^W$ and
 $h_{i,n}\in L^2([0,T])$ satisfying (i) and (ii) of Proposition 2.1 in
 Buckdahn \cite{Bu2}. Fix $n\in{\mathbb N}$ and consider the
 transformation $A^n$ given by (\ref{eq3'a}), when we change $a$ and $A$
 by $a^n$ and $A^n$, respectively. Let $G$ be a smooth functional 
 defined by the right-hand side of (\ref{eq:rsmo}). Then, Buckdahn 
 \cite{Bu} (Proposition 2.2.13) leads to establish \begin{equation}\label{eq4.2mex}
\frac{d}{dt}G(A_{\tau_1^\varepsilon
,t}^n)=-a_t^n D_t^W\left[G(A_{\tau_1^\varepsilon,t}^n)\right].
\end{equation} 
Taking into account (\ref{eqr0mex}), we get
\begin{eqnarray*}
\mathbb{E}_W\left[Y_t^\varepsilon G(A_{\tau_1^\varepsilon,t}^n)\right]
&=&
\mathbb{E}_W\left[
X_{\tau_1^\varepsilon}^\varepsilon 
G(A_{\tau_1^\varepsilon,t}^n)\right]
+\mathbb{E}_W
\left[\int_{\tau_1^\varepsilon}^t a_s Y^\varepsilon_s D_s^W\left[G(A_{\tau_1^\varepsilon,t}^n)\right] ds\right]\\
&&+ 
\mathbb{E}_W\left[\int_{\tau_1^\varepsilon}^t b_s^\varepsilon
 Y^\varepsilon_s G(A_{\tau_1^\varepsilon,t}^n)  ds\right].
\end{eqnarray*}
Now, replacing $G(A_{\tau_1^\varepsilon,t}^n)$ by $G(A_{\tau_1^\varepsilon,s}^n)+\int_{s}^t \frac{d}{dr}
G(A_{\tau_1^\varepsilon,r}^n)dr$,
for some $s\in [\tau_1^\varepsilon, t]$, and using (\ref{eq4.2mex}),
 we obtain
\begin{eqnarray*}
\mathbb{E}_W\left[Y_t^\varepsilon G(A_{\tau_1^\varepsilon,t}^n)\right]
&=&
\mathbb{E}_W\left[X_{\tau_1^\varepsilon}^\varepsilon G\right]+\mathbb{E}_W\left[
X_{\tau_1^\varepsilon}^\varepsilon \int_{\tau_1^\varepsilon}^t \frac{d}{ds}G(A_{\tau_1^\varepsilon,s}^n)ds\right]\\
&&+\mathbb{E}_W\left[
\int_{\tau_1^\varepsilon}^t a_s Y^\varepsilon_s D_s^W\left[G(A_{\tau_1^\varepsilon,s}^n)+\int_{s}^t \frac{d}{dr}G(A_{\tau_1^\varepsilon,r}^n)dr\right] ds\right]\\
&&+ 
\mathbb{E}_W\left[\int_{\tau_1^\varepsilon}^t b_s^\varepsilon
 Y^\varepsilon_s \left(G(A_{\tau_1^\varepsilon,s}^n)
 +\int_{s}^t \frac{d}{dr}G(A_{\tau_1^\varepsilon,r}^n)dr\right)
 ds\right]\\
&=&
\mathbb{E}_W
\left[X_{\tau_1^\varepsilon}^\varepsilon G\right]-\mathbb{E}_W
\left[X_{\tau_1^\varepsilon}^\varepsilon
 \int_{\tau_1^\varepsilon}^t a^n_s D_s^W\left[G(A_{\tau_1^\varepsilon,s}^n)\right]ds\right]\\
&&+ 
\mathbb{E}_W\left[\int_{\tau_1^\varepsilon}^t a_s 
Y^\varepsilon_s D_s^W\left[
G(A_{\tau_1^\varepsilon,s}^n)\right] ds\right]\\
&&-\mathbb{E}_W\left[
\int_{\tau_1^\varepsilon}^t a_s Y^\varepsilon_s 
D_s^W\left[\int_{s}^t a^n_r D_r^W\left(
G(A_{\tau_1^\varepsilon,r}^n)\right)dr\right] 
ds\right]\\
&&+
\mathbb{E}_W\left[
\int_{\tau_1^\varepsilon}^t b_s^\varepsilon Y^\varepsilon_s G(A_{\tau_1^\varepsilon,s}^n)
ds\right]-
\mathbb{E}_W\left[
\int_{\tau_1^\varepsilon}^t b_s^\varepsilon 
Y^\varepsilon_s \int_{s}^t a^n_r
D_r^W\left(G(A_{\tau_1^\varepsilon,r}^n)\right)dr
ds\right].\end{eqnarray*}
Therefore,  the Fubini theorem allows us to state
\begin{eqnarray*}
\mathbb{E}_W\left[Y_t^\varepsilon G(A_{\tau_1^\varepsilon,t}^n)\right]
&=&
\mathbb{E}_W\left[X_{\tau_1^\varepsilon}^\varepsilon G\right]-\mathbb{E}_W\left[
X_{\tau_1^\varepsilon}^\varepsilon \int_{\tau_1^\varepsilon}^t a^n_s D_s^W\left[G(A_{\tau_1^\varepsilon,s}^n)\right]ds\right]\\
&&+ 
\mathbb{E}_W\left[\int_{\tau_1^\varepsilon}^t a_s Y^\varepsilon_s D_s^W\left[G(A_{\tau_1^\varepsilon,s}^n)\right] ds\right]\\
&&-\mathbb{E}_W\left[
\int_{\tau_1^\varepsilon}^t 
\int_{\tau_1^\varepsilon}^r a_s Y^\varepsilon_s D_s^W
\left[ a^n_r D_r^W\left(
G(A_{\tau_1^\varepsilon,r}^n)\right)\right] ds dr\right]\\
&&+
\mathbb{E}_W\left[\int_{\tau_1^\varepsilon}^t b_s^\varepsilon
 Y^\varepsilon_s G(A_{\tau_1^\varepsilon,s}^n)
ds\right]-
\mathbb{E}_W\left[
\int_{\tau_1^\varepsilon}^t \int_{\tau_1^\varepsilon}^r 
b_s^\varepsilon Y^\varepsilon_s a^n_r
D_r^W\left(G(A_{\tau_1^\varepsilon,r}^n)
\right)ds dr\right].
\end{eqnarray*}
By Lemma 2.2.13 in \cite{Bu} we have that $a^n_r D_r^W\left(G(A_{\tau_1^\varepsilon,r}^n)\right)$ 
is a smooth functional for fixed $r\in [\tau_1^\varepsilon,t]$.
 Therefore, applying (\ref{eqr0mex}) to this smooth functional 
 we get that
\begin{eqnarray*}
\mathbb{E}_W\left[Y_t^\varepsilon G(A_{\tau_1^\varepsilon,t}^n)
\right]&=&
\mathbb{E}_W\left[X_{\tau_1^\varepsilon}^\varepsilon G\right]+
\mathbb{E}_W\left[\int_{\tau_1^\varepsilon}^t 
\left(a_s-a_s^n\right) Y_s^\varepsilon D_s^W\left[G(A_{\tau_1^\varepsilon,s}^n)\right]ds\right]\\
&&+\mathbb{E}_W
\left[\int_{\tau_1^\varepsilon}^t b_s^\varepsilon
 Y_s^\varepsilon G(A_{\tau_1^\varepsilon,s}^n)
ds\right].
\end{eqnarray*}
The hypotheses assumed allow us to use the dominated convergence 
theorem. Thus, taking limit as $n\rightarrow \infty$,
\begin{displaymath}
\mathbb{E}\left[Y_t^\varepsilon G(A_{\tau_1^\varepsilon,t})\right]=
\mathbb{E}\left[X_{\tau_1^\varepsilon}^\varepsilon G\right]+
\mathbb{E}\left[\int_{\tau_1^\varepsilon}^t 
b_s^\varepsilon Y_s^\varepsilon G(A_{\tau_1^\varepsilon,s})
ds\right].
\end{displaymath}
From Girsanov's argument  (\ref{egirsanov2}) we have, with
$T_{\tau_1^\varepsilon,s}=A^{-1}_{\tau_1^\varepsilon,s}$,
\begin{displaymath}
\mathbb{E}\left[Y_t^\varepsilon (T_{\tau_1^\varepsilon,t}) 
\mathcal{L}_{\tau_1^\varepsilon,t} G\right]=
\mathbb{E}\left[X_{\tau_1^\varepsilon}^\varepsilon G\right]+
\mathbb{E}\left[\int_{\tau_1^\varepsilon}^t b_s^\varepsilon
(T_{\tau_1^\varepsilon,s})
 Y_s^\varepsilon(T_{\tau_1^\varepsilon,s })  
 \mathcal{L}_{\tau_1^\varepsilon,s} 
G ds\right].\end{displaymath}
Since this is true for any smooth functional $G$, it implies
\begin{displaymath}
Y_t^\varepsilon (T_{\tau_1^\varepsilon,t}) \mathcal{L}_{\tau_1^\varepsilon,t}=
X_{\tau_1^\varepsilon}^\varepsilon
+ \int_{\tau_1^\varepsilon}^t b_s^\varepsilon
(T_{\tau_1^\varepsilon,s})
 Y_s^\varepsilon(T_{\tau_1^\varepsilon,s })  
 \mathcal{L}_{\tau_1^\varepsilon,s} 
ds.
\end{displaymath}
So,
\begin{displaymath}
Y_t^\varepsilon (T_{\tau_1^\varepsilon,t}) \mathcal{L}_{\tau_1^\varepsilon,t}=
X_{\tau_1^\varepsilon}^\varepsilon
\exp\left\{\int_{\tau_1^\varepsilon}^t 
b_s^\varepsilon(T_{\tau_1^\varepsilon,s}) ds\right\},
\end{displaymath}
and, by  (\ref{erelacion1}) and (\ref{eufff}), this means that
\begin{eqnarray*}
Y_t^\varepsilon  &=&X_{\tau_1^\varepsilon}^\varepsilon(
A_{\tau_1^\varepsilon,t})
\exp\left\{\int_{\tau_1^\varepsilon}^t
 b_s^\varepsilon(A_{s,t}) ds\right\} L_{\tau_1^\varepsilon,t}\\
&=&X_0( A_{0,t})
\exp\left\{\int_{0}^{\tau_1^\varepsilon} b_s^\varepsilon
(A_{s,t})\ ds\right\}
\exp\left\{\int_{\tau_1^\varepsilon}^t b_s^\varepsilon
(A_{s,t}) ds\right\}\\
&&\times L_{0,\tau_1^\varepsilon}(A_{\tau_1^\varepsilon,t})
L_{\tau_1^\varepsilon,t}\left[1+v_{\tau_1^\varepsilon-}
(y^\varepsilon_1,A_{\tau_1^\varepsilon,t})\right].
\end{eqnarray*}
Finally, using (\ref{relacion2}), we have that
\begin{displaymath}
Y_t^\varepsilon= X_0( A_{0,t})
\exp\left\{\int_{0}^t b_s^\varepsilon(A_{s,t}) ds\right\} L_{0,t}
\left[1+v_{\tau_1^\varepsilon-}
(y^\varepsilon_1,A_{\tau_1^\varepsilon,t})\right].
\end{displaymath}
This completes the proof for  $t\in [\tau_1^\varepsilon,
\tau_2^\varepsilon)$. The rest of the cases 
can be treated similarly. 
\hfill{$\square$}

\section{Existence and uniqueness of solution for the main equation} 
\setcounter{equation}{0}
\setcounter{theorem}{0}

The main goal of this section is to prove the following theorem.
\begin{theorem}\label{the:snom}
Assume (H1), (H2) and (H4). Suppose also that $a, b, v_{\cdot-}(y)\in L^\infty(\Omega\times[0,T])$, for any $y\in{\mathbb R}_0$, and $D^Wa 
\in L^\infty(\Omega\times[0,T]^2)$.
Then, the process $X$ defined in (\ref{eintr9bis}) is the unique solution  in $L^1(\Omega\times[0,T])$ of 
\begin{equation}\label{enewversion2}
X_t= X_0+
\int_0^t b_s X_s ds+\int_0^t a_s X_s \delta W_s+
 \int_0^t \int_{\mathbb{R}_0} v_{s-}(y) X_{s-}\ d \tilde N(s,y),
 \quad t\in[0,T],
\end{equation}
such that 
$$\left(\int_0^T\int_{{\mathbb R}_0}|v_{s-}(y)X_{s-}|dN(s,y)
+\int_0^T\int_{{\mathbb R}_0}|g(y)X_{s-}|\nu(dy)ds\right)
\in L^1(\Omega).$$
Here, the stochastic integrals with respect to ${\tilde N}$ and
$N$ are pathwise defined.
\end{theorem}

\medskip

\noindent {\bf Remark}. As an immediate consequence  of the proof
of this result,  equation (\ref{eintr2}) can be rewritten as
an  equation of the form (\ref{enewversion2}). 

\medskip

\noindent {\bf Proof of Theorem \ref{the:snom}:} This proof is divided into two parts.

\medskip

\noindent {\bf Step 1.} First of all, by means of a limit argument we will show that $X$ defined in (\ref{eintr9bis}) is  a solution of (\ref{enewversion2}). Towards this end, we 
 prove the convergence of (\ref{emexico7}),
  when $\varepsilon$ tends to zero.

Using that $a$ and $b$ belong to $L^\infty(\Omega\times [0,T])$ 
and that $G$ is a smooth element, we obtain that
\begin{displaymath}
\lim_{\varepsilon \downarrow 0} \mathbb{E}\bigg[\int_0^t a_s X_s^\varepsilon  D_s^WG ds
\bigg]= \mathbb{E}\bigg[\int_0^t a_s X_s  D_s^WG ds
\bigg],
\end{displaymath}
and 
\begin{displaymath}
\lim_{\varepsilon \downarrow 0}\mathbb{E}\bigg[G\Big( X_t^\varepsilon - X_0 -\int_0^t b_s X_s^\varepsilon ds\Big)\bigg]=\mathbb{E}\bigg[G\Big( X_t - X_0 -\int_0^t b_s X_s ds\Big)\bigg].
\end{displaymath}
It only remains to prove that, for any $t\in [0,T]$,
\begin{equation}\label{equnova}
\lim_{\varepsilon \downarrow 0} \mathbb{E}\bigg[G\int_0^t\int_{\{|y|>\varepsilon\}} v_{s-}(y) X_{s-}^\varepsilon\ d \tilde N(s,y)
\bigg]= \mathbb{E}\bigg[G\int_0^t\int_{\{|y|>0\}} v_{s-}(y) X_{s-}\ d \tilde N(s,y)
\bigg].
\end{equation}
In order to prove this convergence 
and that the right-hand side is well-defined,
we utilize the following estimation: 
\begin{displaymath}
\mathbb{E}\left[\left|\int_0^t\int_{\{|y|>\varepsilon\}} v_{s-}(y) X_{s-}^\varepsilon\ d \tilde N(s,y)-
\int_0^t\int_{\{|y|>0\}} v_{s-}(y) X_{s-}\ d \tilde N(s,y)\right|\right]\le I_1^\varepsilon+I_2^\varepsilon,\\
\end{displaymath}
with
\begin{eqnarray*}
I_1^\varepsilon &=&\mathbb{E}\left[ \left|\int_0^t\int_{\{0<|y|\le\varepsilon\}} v_{s-}(y) X_{s-}\ d \tilde N(s,y)\right|\right],\\
I_2^\varepsilon &=&\mathbb{E}\left[\textsc{}\left|\int_0^t\int_{\{|y|>\varepsilon\}} v_{s-}(y) \left[X_{s-}^\varepsilon- X_{s-}\right]\ d \tilde N(s,y)\right|\right].
\end{eqnarray*}
First of all, by the definition of  $\tilde N$, we can write
$$I_1^\varepsilon\le  I_{1,1}^\varepsilon+I_{1,2}^\varepsilon,$$
with
\begin{eqnarray*}
I_{1,1}^\varepsilon &=& \mathbb{E}\left[ \left|\int_0^t\int_{\{0<|y|\le\varepsilon\}} v_{s-}(y) X_{s-}\ d N(s,y)\right|\right],\\
I_{1,2}^\varepsilon &=&\mathbb{E}\left[ \left|\int_0^t\int_{\{0<|y|\le\varepsilon\}} v_{s-}(y) X_{s-}\ \nu(dy) ds\right|\right].
\end{eqnarray*}
Now, (H2) and the bound of $X$ given in the proof of
Proposition \ref{ela1}   imply that 
\begin{eqnarray*}
I_{1,1}^\varepsilon &\le & \mathbb{E}\left[\int_0^T \int_{\{0<|y|\le\varepsilon\}} g(y) \left| X_{s-}\right|\ d N(s,y)\right]
\\
&\le &C\ \mathbb{E}\left[ \int_0^T\int_{\{0<|y|\le\varepsilon\}} g(y)\ L_{0,t}\ \exp\left\{\int_0^T \int_{\mathbb{R}_0} g(y) dN(s,y)\right\} \ d N(s,y)\right]\\
&\le &C\ \mathbb{E}_N\left[ \left(\int_0^T\int_{\{0<|y|\le\varepsilon\}} g(y)\ d N(s,y)\right)\exp\left\{\int_0^T \int_{\mathbb{R}_0} g(y) d N(s,y)\right\}\right].
\end{eqnarray*}
Then, by (H4) we obtain that
\begin{displaymath}\mathbb{E}\left[I_{1,1}^\varepsilon\right]\longrightarrow 0,\end{displaymath}
when $\varepsilon$ goes to 0.

Proceeding similarly, we also get
\begin{displaymath}\mathbb{E}\left[I_{1,2}^\varepsilon\right]\longrightarrow 0,\end{displaymath}
when $\varepsilon$ goes to 0.

Finally, again the relation between of $N$ and $\tilde N$, the fact that
\begin{displaymath}
\left|X_{s-}^\varepsilon -X_{s-}\right|\le 2C\ L_{0,t}\ \exp\left\{\int_0^T \int_{\mathbb{R}_0} g(y) d\ N(s,y)\right\},\end{displaymath}
which is a 
consequence of Proposition \ref{ela1},
and the dominated convergence theorem allow us to ensure that   when $\varepsilon$ goes to 0,
\begin{displaymath}\mathbb{E}\left[I_{2}^\varepsilon\right]\longrightarrow 0.\end{displaymath}
So, the convergence (\ref{equnova}) is satisfied.

\medskip

\noindent {\bf Step 2}. Now we show the uniqueness of the
solution to (\ref{enewversion2}). 
Let $Y\in L^1(\Omega\times[0,T])$ be a
solution of (\ref{enewversion2}) satisfying 
$aY\ \1_{[0,t]}\in Dom\, \delta^W$,
$t\in[0,T]$, and
$$\left(\int_0^T\int_{{\mathbb R}_0}|v_{s-}(y)Y_{s-}|dN(s,y)
\right)
\in L^1(\Omega).$$

 Recall that the coefficients verify that for any $\omega''\in \Omega_N$ a.s.
 \begin{equation}\label{ecoefft}
 \left|b_t(\cdot,w'')\right|+
 \left|a_t(\cdot,w'')\right|+\left|D^W_s a_t(\cdot,w'')\right|+
 \left|v_{t-}(y,\cdot,w'')\right|\le C,
 \end{equation}
for any $s,t \in [0,T]$, $\omega'\in \Omega_W$ and $y\in \mathbb{R}_0$.

Now fix $\omega''\in\Omega_N$, and let $a^n$ and $A^n$ be as in the
proof of Theorem \ref{terss}.
Consequently, for any $G\in \mathcal{S}^W$, we have, by Lemma \ref{lem:tdi},
\begin{eqnarray}\label{equalllv}
\mathbb{E}_W\left[Y_t(\cdot,\omega'') G(A_t^n)\right]&=&
\mathbb{E}_W\left[X_0(\cdot,\omega'') G(A_t^n)\right]+
\mathbb{E}_W\left[G(A_t^n)\int_0^t b_s(\cdot,\omega'') Y_s(\cdot,\omega'') ds\right]\nonumber\\
&&+\mathbb{E}_W\left[\int_0^t a_s(\cdot,\omega'') Y_s(\cdot,\omega'') D^W_s(G(A_t^n)) ds\right]\nonumber\\
&&+ \mathbb{E}_W\left[G(A_t^n)\int_0^t \int_{\mathbb{R}_0} v_{s-}(y,\cdot,\omega'') Y_{s-}(\cdot,\omega'') d\tilde  N(s,y)\right].
\end{eqnarray}
By (2.2.24) in Lemma 2.2.13 of \cite{Bu},
\begin{displaymath}\frac{d}{dt} G(A_t^n)=-a_t^n(\cdot,\omega'') D^W_t(G(A_t^n)),\end{displaymath}
and  it implies that
\begin{displaymath}
G(A_t^n)= G(A_s^n)-\int_s^t a_r^n(\cdot,\omega'') D^W_r(G(A_r^n)) dr.
\end{displaymath}
Taking  this last equality into account, we get 
\begin{displaymath}\begin{array}{l}
 \displaystyle\mathbb{E}_W\left[Y_t(\cdot,\omega'') G(A_t^n)\right] =
\mathbb{E}_W\left[X_0(\cdot,\omega'') G\right]-
\mathbb{E}_W\left[X_0(\cdot,\omega'') \int_0^t a_r^n(\cdot,\omega'') D^W_r(G(A_r^n)) dr\right]\\[4mm]
\qquad \displaystyle +\mathbb{E}_W\left[\int_0^t G(A_s^n) b_s(\cdot,\omega'') Y_s(\cdot,\omega'') ds\right]\\[4mm]
\qquad \displaystyle -\mathbb{E}_W\left[\int_0^t\left(\int_s^t a_r^n(\cdot,\omega'') D^W_r(G(A_r^n))dr\right)
b_s(\cdot,\omega'') Y_s(\cdot,\omega'') ds\right]\\[4mm]
\qquad \displaystyle  +
\mathbb{E}_W\left[\int_0^t a_s(\cdot,\omega'') Y_s(\cdot,\omega'') D^W_s(G(A_s^n)) ds\right]\\[4mm]
\qquad \displaystyle - \mathbb{E}_W\left[\int_0^t D_s^W\left(\int_s^t a_r^n(\cdot,\omega'') D^W_r(G(A_r^n))dr\right)
 a_s(\cdot,\omega'') Y_s(\cdot,\omega'')  ds\right]\\[4mm]
\qquad \displaystyle  +\mathbb{E}_W\left[\int_0^t \int_{\mathbb{R}_0} G(A_s^n) v_{s-}(y,\cdot,\omega'') Y_{s-}(\cdot,\omega'') d\tilde  N(s,y)\right]\\[4mm]
\qquad \displaystyle -
\mathbb{E}_W\left[\int_0^t \int_{\mathbb{R}_0} \left(\int_s^t a_r^n(\cdot,\omega'')
 D^W_r(G(A_r^n))dr\right) v_{s-}(y,\cdot,\omega'') Y_{s-}(\cdot,\omega'') d\tilde  N(s,y)\right].
\end{array}\end{displaymath}

Hence, proceeding as in Step 3 of the proof of Theorem \ref{terss}, we
state
 \begin{eqnarray*}
 \mathbb{E}_W\left[Y_t(\cdot,\omega'') G(A_t^n)\right] &=&
\mathbb{E}_W\left[X_0(\cdot,\omega'') G\right]
 +\mathbb{E}_W\left[\int_0^t G(A_s^n) b_s(\cdot,\omega'') Y_s(\cdot,\omega'') ds\right]\\
 && +\mathbb{E}_W\left[\int_0^t \int_{\mathbb{R}_0} G(A_s^n) v_{s-}(y,\cdot,\omega'') Y_{s-}(\cdot,\omega'') d\tilde  N(s,y)\right]\\
 &&+\mathbb{E}_W\left[\int_0^t \left(a_s-a^n_s\right)(\cdot,\omega'') 
 D^W_s(G(A_s^n)) Y_s(\cdot,\omega'') ds \right].\end{eqnarray*}
 Therefore, proceeding as in the proof of Theorem \ref{terss} again,
 we can write
 
 \begin{eqnarray*}
 \mathbb{E}_W\left[Y_t(\cdot,\omega'') G(A_t)(\cdot,\omega'')\right] &=&
\mathbb{E}_W\left[X_0(\cdot,\omega'') G\right]
 +\mathbb{E}_W\left[\int_0^t G(A_s) b_s(\cdot,\omega'') Y_s(\cdot,\omega'') ds\right]\\
 && +\mathbb{E}_W\left[\int_0^t \int_{\mathbb{R}_0} G(A_s) v_{s-}(y,\cdot,\omega'') Y_{s-}(\cdot,\omega'') d\tilde  N(s,y)\right].\end{eqnarray*}
Thus,  Girsanov Theorem implies
\begin{displaymath}\begin{array}{l}
\displaystyle  \mathbb{E}_W\left[\mathcal{L}_t(\cdot,\omega'') Y_t(T_t(\cdot,\omega''),\omega'') G\right]\\[4mm]
 \displaystyle \qquad =
\mathbb{E}_W\bigg[G\bigg(X_0(\cdot,\omega'') 
 +\int_0^t  b_s(T_s(\cdot,\omega''),\omega'') Y_s(T_s(\cdot,\omega''),\omega'') \mathcal{L}_s(\cdot,\omega'') ds\\[4mm]
 \displaystyle\qquad\quad +\int_0^t \int_{\mathbb{R}_0} v_{s-}(y,T_s(\cdot,\omega''),\omega'') Y_{s-}(T_s(\cdot,\omega''),\omega'') 
 \mathcal{L}_s(\cdot,\omega'')
 d\tilde  N(s,y)\bigg)\bigg],\end{array}\end{displaymath}
 which yields
 \begin{displaymath}\begin{array}{l}
\displaystyle   Y_t(T_t(\omega',\omega''),\omega'') \mathcal{L}_t(\omega',\omega'')\\[3mm]
\displaystyle\qquad  =X_0(\omega',\omega'') 
   +\int_0^t  b_s(T_s(\omega',\omega''),\omega'') Y_s(T_s(\omega,\omega''),\omega'') \mathcal{L}_s(\omega',\omega'') ds\\[3mm]
   \displaystyle\qquad\quad +\int_0^t \int_{\mathbb{R}_0} v_{s-}(y,T_s(\omega',\omega''),\omega'') Y_{s-}(T_s(\omega',\omega''),\omega'') 
 \mathcal{L}_s(\omega',\omega'')
 d\tilde  N(s,y),\qquad \omega'\ {\rm a.s.}
   \end{array}\end{displaymath}
   Since all  factors are mesurable, we can  change $\omega'$ for $\omega''$. So, for almost all $\omega'$,
\begin{equation}\begin{array}{l}
\displaystyle   Y_t(T_t(\omega',\omega''),\omega'') \mathcal{L}_t(\omega',\omega'')\\[3mm]
\displaystyle\qquad  =X_0(\omega',\omega'') 
   +\int_0^t  b_s(T_s(\omega',\omega''),\omega'') Y_s(T_s(\omega,\omega''),\omega'') \mathcal{L}_s(\omega',\omega'') ds\\[3mm]
   \displaystyle\qquad\quad +\int_0^t \int_{\mathbb{R}_0} v_{s-}(y,T_s(\omega',\omega''),\omega'') Y_{s-}(T_s(\omega',\omega''),\omega'') 
 \mathcal{L}_s(\omega',\omega'')
 d\tilde  N(s,y),\qquad \omega''\ {\rm a.s.}\label{e1212dfd}
   \end{array}\end{equation}
  Finally, we only need to observe that 
  either Protter \cite{Pr} (Theorem
  37, page 84), or Bojdecki \cite{Boj} (Theorem 13.12) gives
  \begin{eqnarray*}\lefteqn{
  Y_t(T_t(\omega',\omega''),\omega'') \mathcal{L}_t(\omega',\omega'')}\\
  &=& X_0(\omega)\exp\left\{\int_0^t b_s(T_s(\omega),\omega'')ds
  +\int_0^t\int_{{\mathbb R}_0}v_{s-}(y,T_s(\omega),\omega'')
  d{\tilde N}(s,y)\right\}\\
  &&\times\prod_{0\le s\le t}\left[1+v_{s-}(y,T_s(\omega),\omega'')
  \Delta N(s,y)\right],\end{eqnarray*}
  wich means that $X=Y$ and the proof is complete.
\hfill{$\square$}

\section{Appendix}\label{sec:ape}
\setcounter{equation}{0}
\setcounter{theorem}{0}

 This section is devoted to present the proof of Theorem \ref{continuity}. In order to simplify the notation, we use the convention $D=D^W$ because, in this section, the probability space is the
canonical Wiener space. 
 Also, remember that $c$ will denote a generic constant that may change from line to line. 

We begin this section with an auxiliary result.

\begin{lemma}\label{lem:aux}
Under the assumptions of Theorem \ref{continuity}, we have that, for $s\in[0,T]$,
\begin{itemize}
\item[(a)]
$$\int_0^T |D_{\theta}(a_r(A_{r,s}))|^2 d\theta\leq 2e^{2c_1}||\,|Da_r|_2\,||_{\infty}^2,\quad r\in[0,s],$$
where $c_1:=\int_0^T ||\, |Da_r|_2^2\, ||_{\infty} dr$.  
\item[(b)]
$$\int_0^s \int_0^T |D_{\theta}(a_r(A_{r,s}))|^2 d\theta dr\leq 2c_1e^{2c_1}.$$
\item[(c)]
$$\int_0^T |(D_{\theta}a_r)(A_{r,t})-(D_{\theta}a_r)(A_{r,s})|^2d\theta
\leq|| |DD a_r|_2^2 ||_{\infty} 2 e^{2c_1}\int_s^t ||a_r||^2_{\infty} dr
,\quad r\in[0,s].$$
\item[(d)]
$$\int_0^s\int_0^T |(D_{\theta}a_r)(A_{r,t})-(D_{\theta}a_r)(A_{r,s})|^2d\theta dr\leq 2 e^{2c_1}c_2\int_s^t ||a_r||^2_{\infty} dr,$$
with $c_2:=\int_0^T ||\, |DDa_r|_2^2\, ||_{\infty} dr.$
\end{itemize}
\end{lemma}
\noindent {\bf Proof}:
We first observe that, by virtue of Proposition \ref{propbeta} and 
(\ref{eq3'a}), we obtain
\begin{equation}\label{chainrule}
D_{\theta}(a_r(A_{r,s}))=(D_{\theta} a_r)(A_{r,s})-\int_r^s (D_u a_r)(A_{r,s}) D_{\theta}(a_u(A_{u,s}))du.
\end{equation}
Therefore, from Hölder inequality, we have 
\begin{eqnarray*}
&&\int_0^T |D_{\theta}(a_r(A_{r,s}))|^2 d\theta\\
&&\qquad \leq 2\int_0^T |(D_{\theta}a_r)(A_{r,s})|^2 d\theta + 2\int_0^T \left|\int_r^s (D_u a_r)(A_{r,s}) D_{\theta}(a_u(A_{u,s}))du\right|^2 d\theta\\
&&\qquad \leq 2||\,|Da_r|_2\,||_{\infty}^2+2\int_0^T \left(\int_r^s |(D_u a_r)(A_{r,s})|^2 du\right)\left(\int_r^s |D_{\theta}a_u(A_{u,s})|^2 du\right)d\theta\\
&& \qquad \leq 2||\,|Da_r|_2\,||_{\infty}^2+2||\,|Da_r|_2\,||_{\infty}^2 \int_r^s \left(\int_0^T |D_{\theta}(a_u(A_{u,s}))|^2d\theta\right)du.
\end{eqnarray*}
Consequently, by Gronwall's lemma, we deduce
$$\int_0^T |D_{\theta}(a_r(A_{r,s}))|^2 d\theta\leq 2||\,|Da_r|_2\,||_{\infty}^2 \exp\left\{\int_r^s 2||\,|Da_{u}
|_2\,||_{\infty}^2 du\right\}\leq 2e^{2c_1}||\,|Da_r|_2\,||_{\infty}^2,$$
which  shows that Statement (a) is satisfied. 

Now, using Proposition \ref{firstinequality}, Lemma \ref{uniformcm} and the definition of constant $c_1$ we obtain 
\begin{eqnarray*}
\int_0^T |(D_{\theta}a_r)(A_{r,t})-(D_{\theta}a_r)(A_{r,s})|^2d\theta
&\leq & || |DDa_r|_2^2 ||_{\infty} \sup_{r\le s} |A_{r,t}-A_{r,s}|^2_{CM}\\
&\leq & || |DD a_r|_2^2 ||_{\infty} 2 e^{2c_1}\int_s^t ||a_u||^2_{\infty} du.
\end{eqnarray*}
Thus, Statement  (c) holds. 

Finally, Statements (b) and  (d) are an immediate consequence of (a) and (b), and the proof is complete. 
\hfill{$\square$}
\vspace{0.5cm}

Now, we are ready to prove Theorem \ref{continuity}
\vspace{0.5cm}

\noindent {\bf Proof of Theorem \ref{continuity}:} From (\ref{edensitya}) and (\ref{eintr8}), we only need to show the continuity  of processes $Z_{0}(A_{0,\cdot}),$ 
$\int_{0}^\cdot h_s (A_{s,\cdot})ds,$ $\int_{0}^\cdot \ a_s^2(A_{s,\cdot})\ ds,$ $\int_{0}^\cdot a_s(A_{s,\cdot})\ \delta W_s,$ and $\int_{0}^\cdot\int_{s}^\cdot (D_u a_s)(A_{s,\cdot})
\ D_s(a_u(A_{u,\cdot}))\ du ds.$ So now we divide the proof in several steps and we assume that $0\le s\le t\le T$.
\begin{itemize}
\item[{\bf (1)}] Taking into account Proposition \ref{firstinequality}, we get
$$|Z_{0}(A_{0,t})-Z_{0}(A_{0,s})|\leq \left\|\left(\int_{0}^T |D_sZ_{0}|^2 ds\right)^{\frac{1}{2}}\right\|_{\infty} |A_{0,t}\omega-A_{0,s}\omega|_{CM},$$
which, together with the definition of the space $\mathbb{D}^W_{1,\infty}$ and Lemma \ref{uniformcm}, implies that
the process $\{Z_{0}(A_{0,t}):t\in[0,T]\}$ has continuous paths.
\item[{\bf (2)}] We show now the continuity of $\{\int_{0}^t g_r (A_{r,t})dr:t\in[0,T]\}$, with $g\in L^1([0,T], \mathbb{D}_{1,\infty}).$ Note that, in this proof, $g_r:=h_r$ or $g_r:=a_r^2.$ 
  From Proposition \ref{firstinequality}, we can conclude  
\begin{eqnarray*}\lefteqn{
\left|\int_{0}^t g_r (A_{r,t})dr-\int_{0}^s g_r (A_{r,s})dr\right|}\\
&&\qquad \leq \int_{s}^t ||g_r||_{\infty}dr+\int_{0}^s ||\,|Dg_r|_2\,||_{\infty} |A_{r,t}-A_{r,s}|_{CM} dr\\
&&\qquad \leq\int_{s}^t ||g_r||_{\infty}dr
+\left(\int_{0}^T ||\,|Dg_r|_2\,||_{\infty} dr\right) \sup_{r\le s} |A_{r,t}-A_{r,s}|_{CM},
\end{eqnarray*}
which gives the desired continuity due to Lemma \ref{uniformcm}. 
\item[{\bf (3)}] Next step is to check the continuity of 
$\{\int_{0}^t \ a_r(A_{r,t})\ \delta W_r:t\in [0,T]\}$. So, 
 by the Kolmogorov-Centsov continuitiy criterion (see \cite{N06}, for example), it is sufficient to show that for some $p\in (2,\infty),$ 
$${\mathbb E}\left|\int_{0}^t \ a_r(A_{r,t})\ \delta W_r-\int_{0}^s \ a_r(A_{r,s})\ \delta W_r\right|^{2p}\leq c (t-s)^{p-1},$$
where $c$ is a constant that only depends on $p$ and $a$.

Remember that $a_\cdot(A_{\cdot,t}){1\!\!1}_{[0,t]}(\cdot)$ belongs to $L^2 ([0,T],{\mathbb D}_{1,2})$ as a consequence of Propositions \ref{propbeta} and \ref{propgamma} (see also \cite{Bu,Bu2,bu3}). In particular this guarantee that $\delta(a_\cdot (A_{\cdot,t}){1\!\!1}_{[0,t]}(\cdot))$ is well-defined. 
We can apply  H\"older inequality and Proposition 3.2.1 in \cite{N06} to derive
\begin{eqnarray*}\lefteqn{
{\mathbb E}\left|\int_{0}^t \ a_r(A_{r,t})\ \delta W_r-\int_{0}^s \ a_r(A_{r,s})\ \delta W_r\right|^{2p}}\\
&\leq& c\ {\mathbb E}\left(\left|\int_s^t \ a_r(A_{r,t})\ \delta W_r\right|^{2p}\right)+c\ {\mathbb E}\left(\left|\int_{0}^s \ [a_r (A_{r,t})-a_r(A_{r,s})]\ \delta W_r\right|^{2p}\right)\\
&\leq& c \left(\int_s^t ({\mathbb E}(a_r(A_{r,t})))^2dr\right)^p+c\
{\mathbb E}\left(\int_s^t \int_0^T (D_{\theta} a_r (A_{r,t}))^2 d\theta dr\right)^p\\
&& +c \left(\int_{0}^s ({\mathbb E}(a_r(A_{r,t})-a_r(A_{r,s})))^2dr\right)^p\\
&&+
c\ {\mathbb E}\left(\int_{0}^s \int_0^T (D_{\theta} (a_r (A_{r,t})-a_r (A_{r,s})))^2 d\theta dr\right)^p\\
&=&c\ \{A+B+C+D\}.
\end{eqnarray*}

In order to finish this step, we are going to see that these four terms are bounded by $c (t-s)^{p-1}.$ Towards this end, 
we observe that Hölder inequality 
and Lemma \ref{lem:aux} (a)
allow us to conclude 
$$A\leq \left(\int_0^T ({\mathbb E}(a_r(A_{r,t})))^{2p} dr\right)(t-s)^{p-1}\leq |\,||a_r||_{\infty}\,|_{2p}^{2p}\, (t-s)^{p-1}$$
and 
$$B\leq c\left(\int_0^T ||\,|Da_r|_2\,||_{\infty}^{2p} dr\right)(t-s)^{p-1}.$$
Using Hölder inequality again, together with Proposition \ref{firstinequality}, Lemma \ref{uniformcm} and the definition of $c_1$ given in the previous lemma, we have
\begin{eqnarray*}
C
&\leq&\left(\int_{0}^s {\mathbb E}(|a_r(A_{r,t})-a_r(A_{r,s})|^2)dr\right)^p\\
&\leq&\left(\int_{0}^s ||\,|Da_r|_2\,||_{\infty}^2 {\mathbb E}(|A_{r,t}-A_{r,s}|^2_{CM})dr\right)^p\\
&\leq&\left(\int_{0}^s ||\,|Da_r|_2\,||_{\infty}^2 2\left(\int_s^t ||a_u||^2_{\infty} du\right)\exp\left\{2\int_r^s ||\, |Da_u|_2^2\, ||_{\infty} du\right\}  dr\right)^p\\
&\leq&2^p c_1^p e^{2pc_1}\left(\int_s^t ||a_\theta||^2_{\infty} d\theta\right)^p\\
&\leq&2^p c_1^p e^{2pc_1}\left(\int_{0}^T ||a_\theta||^{2p}_{\infty} d\theta\right) (t-s)^{p-1}.
\end{eqnarray*}

In order to manage term $D$, we observe that equation (\ref{chainrule}), Lemma \ref{lem:aux} and Cauchy-Schwarz inequality lead to establish
\begin{eqnarray*}
&&\int_0^T |D_{\theta}[a_r(A_{r,t})-a_r(A_{r,s})]|^2 d\theta\\
&&\quad \leq 2\int_0^T |(D_{\theta}a_r)(A_{r,t})-(D_{\theta}a_r)(A_{r,s})|^2d\theta\\
&&\qquad+ 4\left(\int_s^t |(D_ua_r)(A_{r,t})|^2 du\right)\int_s^t\int_0^T |D_{\theta}(a_u (A_{u,t}))|^2 d\theta du\\
&&\qquad+8\left(\int_r^s \int_0^T |D_{\theta}(a_u(A_{u,t}))|^2 d\theta du\right)\int_r^s |(D_u a_r)(A_{r,t})-(D_u a_r)(A_{r,s})|^2 du\\
&&\qquad+ 8\left(\int_r^s |(D_ua_r)(A_{r,s})|^2 du\right)\int_r^s \int_0^T |D_{\theta} (a_u (A_{u,t}))-D_{\theta}(a_u (A_{u,s}))|^2 d\theta du\\
&&\quad \leq4 e^{2c_1} || |DDa_r|_2^2 ||_{\infty}\int_s^t ||a_u||^2_{\infty} du
+8 e^{2c_1} || |Da_r|_2^2 ||_{\infty} \int_s^t || |Da_u|_2^2 ||_{\infty} du\\
&&\qquad+32 e^{4c_1} c_1 || |DDa_r|_2^2 ||_{\infty} \int_s^t ||a_u||^2_{\infty} du\\
&&\qquad+ 8|| |Da_r|_2^2 ||_{\infty} \int_r^s \int_0^T |D_{\theta} (a_u (A_{u,s}))-D_{\theta}(a_u (A_{u,t}))|^2 d\theta du.
\end{eqnarray*}
Let $v\in [0,s].$ Joining first and third terms in the right hand side of the last expression and integrating in both sides with respect to $r$ between $v$ and $s$ we obtain
\begin{eqnarray*}
&&\int_v^s \int_0^T |D_{\theta}[a_r(A_{r,t})-a_r(A_{r,s})]|^2 d\theta dr\\
&&\quad \leq {\tilde c}_1 \left(\int_v^s || |DDa_r|_2^2 ||_{\infty} dr\right)\int_s^t ||a_r||^2_{\infty} dr\\
&&\qquad+8 e^{2c_1} \left(\int_v^s || |Da_r|_2^2 ||_{\infty}dr\right)\left(\int_s^t || |Da_r|_2^2 ||_{\infty}dr\right)\\
&&\qquad+8\int_v^s || |Da_r|_2^2 ||_{\infty}\int_r^s \int_0^T |D_{\theta}[a_u(A_{u,t})-a_u(A_{u,s})]|^2 d\theta du dr,
\end{eqnarray*}
where ${\tilde c}_1:=4 e^{2c_1}+32 e^{4c_1} c_1.$ 

Now, applying Gronwall's lemma, 
\begin{equation}\label{fitaD}
\int_v^s \int_0^T (D_{\theta} (a_r (A_{r,t})-a_r (A_{r,s})))^2 d\theta dr \leq c \left\{\int_s^t ||a_r||_{\infty}^2 dr+\int_s^t || |Da_r|_2^2 ||_{\infty} dr\right\}.
\end{equation}
Therefore, using Minkowski and Hölder inequalities we state 
$$
D \leq 2^{p-1} c \left(\int_0^T ||a_r||^{2p}_{\infty} dr+\int_0^T ||
 |Da_r|_2 ||_{\infty}^{2p} dr\right)(t-s)^{p-1}.
$$
Thus, the claim of this step  is satisfied.
\item[{\bf (4)}] Finally, we consider the process
 $t\mapsto \int_{0}^t\int_{s}^t (D_u a_s)(A_{s,t})\ D_s(a_u(A_{ut}))\ du ds.$
 
We have 
\begin{eqnarray*}
&&\left|\int_{0}^t\int_{r}^t (D_u a_r)(A_{r,t})\ D_r(a_u(A_{ut}))\ du dr-\int_{0}^s\int_{r}^s (D_u a_r)(A_{r,s})\ D_r(a_u(A_{u,s}))\ du dr\right|\\
&&\quad \leq\int_{s}^t\int_{r}^t |(D_u a_r)(A_{r,t})\ D_r(a_u(A_{u,t}))|du dr\\
&&\qquad+\int_{0}^s\int_{s}^t |(D_u a_r)(A_{r,t})\ D_r(a_u(A_{u,t}))|du dr\\
&&\qquad+\int_{0}^s\int_{r}^s |[(D_u a_r)(A_{r,t})-(D_u a_r)(A_{r,s})] D_r(a_u(A_{u,t}))|du dr\\
&&\qquad+\int_{0}^s\int_{r}^s |(D_u a_r)(A_{r,s})[D_r(a_u(A_{u,t}))-D_r(a_u(A_{u,s}))]|du dr\\
&&\quad =J_1+J_2+J_3+J_4.
\end{eqnarray*}
Now, Cauchy-Schwarz inequality and Lemma \ref{lem:aux} yield 
\begin{eqnarray*}
J_1&\leq&\left(\int_{s}^t\int_{0}^T |(D_u a_r)(A_{r,t})|^2 du dr\right)^{\frac{1}{2}} \left(\int_{0}^T\int_{0}^t |(D_r(a_u(A_{u,t}))|^2 du dr\right)^{\frac{1}{2}}\\
&\leq&\sqrt{2c_1e^{2c_1}}\left(\int_{s}^t ||\, |Da_r|_2^2\, ||_{\infty} dr\right)^{\frac{1}{2}}, \\
J_2&\leq&\left(\int_{0}^s\int_{s}^t |(D_u a_r)(A_{r,t})|^2du dr\right)^{\frac{1}{2}} \left(\int_{0}^s\int_{s}^t |D_r(a_u(A_{u,t}))|^2 du dr\right)^{\frac{1}{2}}\\
&\leq&c_1\sqrt{2e^{2c_1}}\left(\int_{s}^t ||\,|Da_u|_2^2\,||_{\infty} du\right)^{\frac{1}{2}}
\end{eqnarray*}
and 
\begin{eqnarray*}
J_3&\leq&\left(\int_{0}^s\int_{0}^s |(D_u a_r)(A_{r,t})-(D_u a_r)(A_{r,s})|^2 du dr\right)^{\frac{1}{2}}
\left(\int_{0}^s\int_{0}^s |D_r(a_u(A_{u,t}))|^2 du dr\right)^{\frac{1}{2}}\\
&\leq&2\sqrt{c_1c_2} e^{2c_1} \left(\int_s^t ||a_r||_{\infty}^2 dr\right)^{\frac{1}{2}}.
\end{eqnarray*}
Finally, by means of  inequaltiy (\ref{fitaD}), we have
\begin{eqnarray*}
J_4&\leq& \left(\int_{0}^s\int_{r}^s |D_r(a_u(A_{u,t}))-D_r(a_u(A_{u,s}))|^2 du dr\right)^{\frac{1}{2}} 
\left(\int_{0}^s\int_{r}^s |(D_u a_r)(A_{r,s})|^2du dr\right)^{\frac{1}{2}}\\
&\leq& \sqrt{c_1}\left(\int_{0}^s\int_{0}^T |D_r(a_u(A_{u,t}))-D_r(a_u(A_{u,s}))|^2 dr du\right)^{\frac{1}{2}}\\
&\leq& c \left(\int_s^t ||a_r||_{\infty}^2 dr+\int_{s}^t ||\,|Da_r|_2^2\,||_{\infty} dr\right)^{\frac{1}{2}}.
\end{eqnarray*}
Thus, the proof is complete.
\hfill{$\square$}

\end{itemize}

\textbf{Aknowledgements}
The authors thank Cinvestav-IPN and Universitat de Barcelona for their 
hospitality and economical support. The paper was partially supported by the CONACyT grant 98998, MEC FEDER MTM 2009-07203 and
MEC FEDER MTM 2009-08869.

\end{document}